\documentstyle[amssymb,11pt]{article}
\def\G{\Gamma}
\def\D{\Delta}

\def\gO{\Omega}

\def\a{\alpha}
\def\b{\beta}
\def\d{\delta}
\def\k{\kappa}
\def\e{\epsilon}
\def\f{\varphi}
\def\g{\gamma}

\def\h{\eta}
\def\l{\lambda}
\def\m{\mu}
\def\n{\nu}
\def\go{\omega}
\def\p{\pi}

\def\t{\tau}
\def\r{\rho}
\def\s{\sigma}
\def\th{\theta}
\def\x{\xi}
\def\z{\zeta}

\def\bD{{\bf D}}

\def\bG{{\bf G}}
\def\bL{{\bf L}}

\def\bR{{\bf R}}
\def\bC{{\bf C}}
\def\bH{{\bf H}}
\def\bK{{\bf K}}
\def\bN{{\bf N}}

\def\bV{{\bf V}}

\def\bc{{\bf c}}

\def\pf{\rightarrow}

\font\tengoth=eufm10 at 12pt
\font\sevengoth=eufm7 at 8pt
\newfam\gothfam
\textfont\gothfam=\tengoth
\scriptfont\gothfam=\sevengoth

\def\frak{\fam\gothfam\tengoth}

\def\fz{{\frak z}}

\def\fa{{\frak a}}

\def\fe{{\frak e}}

\def\fg{{\frak g}}
\def\fh{{\frak h}}

\def\fk{{\frak k}}
\def\fl{{\frak l}}
\def\fm{{\frak m}}
\def\fn{{\frak n}}
\def\fo{{\frak o}}
\def\fq{{\frak q}}
\def\fp{{\frak p}}

\def\fs{{\frak s}}

\def\fu{{\frak u}}
\def\fz{{\frak z}}

\def\bbC{{\Bbb C}}

\def\bbR{{\Bbb R}}

\def\cA{\mathord{{\cal A}}}

\def\cC{\mathord{{\cal C}}}
\def\cD{\mathord{{\cal D}}}

\def\cH{\mathord{{\cal H}}}

\def\cL{\mathord{{\cal L}}}

\def\cO{\mathord{{\cal O}}}

\def\GL{\mathop{{\rm GL}}\nolimits}
\def\lSL{\mathop{{\rm SL}}\nolimits}

\def\SO{\mathop{{\rm SO}}\nolimits}
\def\SU{\mathop{{\rm SU}}\nolimits}

\def\Sp{\mathop{{\rm Sp}}\nolimits}

\def\lsl{\mathop{{\fs\fl}}\nolimits}
\def\sp{\mathop{{\fs\fp}}\nolimits}
\def\su{\mathop{{\fs\fu}}\nolimits}
\def\so{\mathop{{\fs\fo}}\nolimits}

\def\ad{\mathop{{\rm ad}}\nolimits}
\def\Ad{\mathop{{\rm Ad}}\nolimits}

\def\1/2{\mathop{\frac{1}{2}}\nolimits}
\def\Re{\mathop{{\rm Re}}\nolimits}
\def\GH{\mathop{G\bftdiagup H}\nolimits}

\def\pr{\mathop{{\rm pr}}\nolimits}

\def\Im{\mathop{{\rm Im}}\nolimits}

\def\Hom{\mathop{{\rm Hom}}\nolimits}

\def\Aut{\mathop{{\rm Aut}}\nolimits}

\def\id{\mathop{{\rm id}}\nolimits}

\def\spec{\mathop{{\rm spec}}\nolimits}
\def\Supp{\mathop{{\rm Supp}}\nolimits}

\def\Tr{\mathop{{\rm Tr}}\nolimits}

\def\conv{\mathop{{\rm conv}}\nolimits}

\newcommand{\npar}{\par\hspace*{1.1em}}

\def\*{\times}
\def\+{\oplus}
\def\und{\mathop{\quad\mbox{and}\quad}\nolimits}

\def\off{\offinterlineskip}

\def\ecke{\hbox{$\bullet$}}
\def\cecke{\hbox{$\circ$}}
\def\oecke#1{\off\vbox{\hbox{$#1$}\vskip 1mm\ecke}}
\def\uecke#1{\off\vbox{\ecke\hbox{$#1$}}}

\def\coecke#1{\off\vbox{\hbox{$#1$}\vskip 1mm\cecke}}

\def\senk{\moveright 2.5pt\hbox{\vrule width 1.5pt height 1truecm depth0pt}}

\def\graphs#1,#2{{\off\vtop{\oecke{#1}\senk\uecke{#2} }}}
\def\cgraphs#1,#2{{\off\vtop{\coecke{#1}\senk\uecke{#2} }}}

\DeclareMathSymbol{\bftdiagup}       {\mathord}{AMSb}{"1E}
\DeclareMathSymbol{\bftdiagdown}     {\mathord}{AMSb}{"1F}
\newbox\ipbox
\newcommand{\ip}[2]{\left\langle #1\mathrel{\mathchoice
{\setbox\ipbox=\hbox{$\displaystyle \mathstrut #1#2$}
\vrule height\ht\ipbox width0.25pt depth\dp\ipbox}
{\setbox\ipbox=\hbox{$\textstyle \mathstrut #1#2$}
\vrule height\ht\ipbox width0.25pt depth\dp\ipbox}
{\setbox\ipbox=\hbox{$\scriptstyle \mathstrut #1#2$}
\vrule height\ht\ipbox width0.25pt depth\dp\ipbox}
{\setbox\ipbox=\hbox{$\scriptscriptstyle \mathstrut #1#2$}
\vrule height\ht\ipbox width0.25pt depth\dp\ipbox}
} #2\right\rangle}
\newcommand{\halo}[1]{\bar{#1}}
\newlength{\customskipamount}
\newlength{\customleftmargin}
\newtheorem{Def}{Definition}[section]
\newtheorem{Satz}[Def]{Theorem}
\newtheorem{Lemma}[Def]{Lemma}
\newtheorem{Cor}[Def]{Corollary}
\newtheorem{Rem}[Def]{Remark}
\newtheorem{Ex}[Def]{Example}
\newtheorem{Prop}[Def]{Proposition}
\begin{document}
\title{Unitary Representations of Lie Groups with
Reflection Symmetry}
\author{Palle E. T. Jorgensen\thanks{Work supported in part by 
the U.S. National Science Foundation.}\\
Department of Mathematics\\
University of Iowa\\
Iowa City\\
IA 52242\\
\makebox[0pt]{\hss\small E-mail: \ttfamily jorgen@math.uiowa.edu\hss}
\and
Gestur \'Olafsson\thanks{The author was supported by NSF
grant DMS-9626541, LEQSF grant (1996-99)-RD-A-12 and
the Danish Research Council}
\\Department of Mathematics\\
Louisiana State University\\
Baton Rouge\\
LA 70803\\
\makebox[0pt]{\hss\small E-mail: \ttfamily olafsson@math.lsu.edu\hss}
}
\date{}
\maketitle

\begin{abstract}
We consider the following class of unitary representations $\p $ of
some (real) Lie group $G$ which has a matched pair of symmetries
described as follows: (i) Suppose $G$ has a period-$2$ automorphism
$\t $, and that the Hilbert space $\bH (\p )$ carries a unitary
operator $J$ such that $J\p =(\p \circ \t )J$ (i.e., {\em
selfsimilarity\/}). (ii) An added symmetry is implied if $\bH (\p )$
further contains a closed subspace $\bK_0 $ having a certain {\em
order-covariance} property, and satisfying the $\bK_0 $-restricted {\em
positivity\/}: $\ip{v}{Jv}\ge 0$, $\forall v\in \bK_0 $, where
$\ip{\cdot }{\cdot }$ is the inner product in $\bH (\p )$. {}From
(i)--(ii), we get an induced dual representation of an associated dual
group $G^c$. All three properties, selfsimilarity, order-covariance,
and positivity, are satisfied in a natural context when $G$ is
semisimple and hermitean; but when $G$ is the $(ax+b)$-group, or the Heisenberg
group, positivity is incompatible with the other two axioms for the
infinite-dimensional irreducible representations. We describe a class
of $G$, containing the latter two, which admits a classification of
the possible spaces $\bK_0 \subset \bH (\p )$ satisfying the axioms of
selfsimilarity and order-covariance.
%
\end{abstract}

\section{\protect\label{Introduction}Introduction}
\setcounter{equation}{0}

We consider a class of unitary representations of 
a Lie group $G$ which possess a certain reflection symmetry
defined as follows: If $\p$ is a representation of $G$
in some Hilbert space $\bH$, we introduce the following three structures:
\begin{list}{}{\setlength{\leftmargin}{\customleftmargin} 
\setlength{\itemsep}{0.5\customskipamount}
\setlength{\parsep}{0.5\customskipamount}
\setlength{\topsep}{\customskipamount}}
\item[\hss\llap{\rm i)}]  $\t \in \Aut (G)$ of period $2$;

\item[\hss\llap{\rm ii)}]  $J :\bH\pf \bH$ is a unitary operator of period $2$
such that $J\p (g)J^* =
\p (\t (g))$, $g\in G$
(this will hold if $\p$ is of the form $\p_+\oplus \p_-$
with $\p_+$ and $\p_-\circ \t$ unitarily equivalent); it
will further be assumed that there is a closed subspace
$\bK_0\subset \bH$ which is invariant under 
$\p (H)$, $H = G^\t$, or more generally, an open subgroup of
$G^\t$;

\item[\hss\llap{\rm iii)}]  positivity is assumed in the sense that
$\ip{v}{J(v)}\ge 0$, $v\in \bK_0$; and
%
\end{list}

\npar
Let $\fg$ be
the Lie algebra of $G$, and let $\fh$ be the Lie algebra
of the fixed-point subgroup $G^\t = \{g\in G\mid \t (g) =g\}$. Let $\fq =
\{Y \in \fg \mid \t(Y) = -Y\}$. Then
\[ \fg =\fh \oplus \fq\, .\]
Let $H$ be a closed subgroup of $G$, $G^\t_o\subset H
\subset G^\t$. 
Assume there
is an $H$-invariant, closed, and generating convex cone $C$ in $\fq$
(i.e., $C-C = \fq$) such that $C^o$ consists of hyperbolic elements.
We assume that $S (C) = H\exp C$ is a closed semigroup in $G$
which is
homeomorphic to $H\times C$, and that
\[H\times C^o\ni (h,Y)\mapsto h\exp Y \in S^o\]
is a diffeomorphism.

\npar We shall consider closed subspaces $\bK_0\subset \bH(\p)$, where
$\bH (\p)$ is the Hilbert space of $\p$ such that $\bK_0$ is invariant
under $\p (S^o)$. Let $J : \bH (\p)\pf \bH (\p)$ be a unitary
period-two intertwiner for $\p$ and $\p \circ \t$, and assume that
$\bK_0$ may be chosen such that $\|v\|^2_J := \ip{v}{Jv} \ge 0$ for
all $v\in \bK_0$.  We will always assume our inner product conjugate
linear in the first argument. We form, in the usual way, the Hilbert
space $\bK = \left(\bK_0\bftdiagup \bN\right)\tilde{}$ by dividing out
with $\bN = \{v\in \bK_0\mid \ip{v}{Jv} = 0\}$ and completing in the
norm $\| \cdot \|_J$.  (This is of course a variation of the
Gelfand-Naimark-Segal (GNS) construction.)  With the properties of
$(G,\p , \bH (\p),\bK)$ as stated, we show using the L\"{u}scher-Mack
theorem that the simply connected Lie group $G^c$ with Lie algebra
$\fg^c = \fh \oplus i\fq$ carries a {\it unitary} representation
$\p^c$ on $\bK$ such that $\{\p^c(h\exp (iY))\mid h\in H, Y\in C^o\}$
is obtained from $\p$ by passing the corresponding operators $\p
(h\exp Y)$ to the quotient $\bK_0\bftdiagup \bN$. In fact, when $Y\in
C$, the selfadjoint operator $d\p (Y)$ on $\bK$ has spectrum contained
in $(-\infty , 0]$.  As in Corollary \ref{MNoughtZero}, we show that
in the case where $C$ extends to an $G^c$ invariant regular cone in
$i\fg^c=i\fh\oplus \fq$ and $\p^c$ is injective, then each $\p^c$ (as
a unitary representation of $G^c$) must be a direct integral of
highest-weight representations of $G^c$.  The examples show that one
can relax the condition in different ways, i.e., one can avoid using
the L\"{u}scher-Mack theorem by instead constructing local
representations and using only cones that are neither generating nor
$H$-invariant.

\npar
Let us outline the plan by a simple examples. Let $G =
\lSL(2,\bbR)$, and let $P$ be the parabolic subgroup
\[P =\left\{\left. p(a,x)=\left(\matrix{ a & x\cr
0 & a^{-1}\cr}\right)\,\right|\,a\in \bbR^*,\, x\in \bbR\right\}\, 
.\]
For $s\in \bbC$, let $\p_s$ be the representation of
$G$ acting by $[\p_s(a)f](b) = f(a^{-1}b)$ on the space
$\bH_s$
of functions $f : G\pf \bbC$,
\[f(gp(a,x) ) = |a|^{-s-1}f(g)\, ,\quad\int_{\SO(2)}|f(k)|^2\,dk<\infty\, ,\]
and with inner product
\[\ip{f}{g} = \int_{\SO(2)}\overline{f(k)}g(k)\,dk\, ,\]
i.e., $\p_s$ is the principal series representation of
$G$ with parameter $s$. The representations $\p_s$ are
unitary in the above Hilbert-space structure as long
as $s\in i\bbR$. For defining a unitary structure for
other parameters we need the intertwining
operator $A_s : \bH_s \pf \bH_{-s}$ defined by
\[A_s(f)(g) = \int_{-\infty}^{\infty}f (gw\halo{n}_y)\,dy\]
for $\Re s\ge 0$ and then generally by analytic continuation. Here $w$
is the Weyl group element $w = \left( \matrix{ 0 & 1\cr -1 &
0\cr}\right)$ and $\halo{n}_y =\left( \matrix{ 1 & 0\cr y &
1\cr}\right)$.

\npar
By restriction to $\halo{N} = \{\halo{n}_y\mid y\in 
\bbR\}$ we can also realize the representations $\p_s$ on
$\bbR\simeq \halo{N}, \, y\mapsto \halo{n}_y$. Using
that
\[\left(\matrix{ 1 & 0\cr y & 1\cr}\right) 
\left(\matrix{ \a & 0\cr 0 & \a^{-1}\cr}\right) 
\left(\matrix{ 1 & x\cr 0 & 1\cr}\right)
=  \left(\matrix{ \a & \a x\cr\a y & \a yx + 
\a^{-1}\cr}\right)\]
we get that $g = \left(\matrix{ a & b\cr c & d\cr}\right)
\in \halo{N}P$ if and only if $a\not= 0$, and in that case
\begin{equation}\label{NbarPmax}
\left(\matrix{ a & b\cr c & d\cr}\right)
= \halo{n}_{c/a}p(a,b/a)\, .
\end{equation}
Thus the intertwiner $A_s$ becomes the singular integral operator
\[A_sf(x) = \int_{-\infty}^\infty f(y)|x - 
y|^{s-1} \,dy\, .\]
In the new inner product
\[ \ip{f}{A_sg} = 
\int_{-\infty}^{\infty}\int_{-\infty}^{\infty}\overline{f(x)}g(y)|x-y|^{s-1}\, 
dx\,dy\]
the representation $\p_s$, which is given by

\[\left[\p_s\left(\left(\matrix{ a & b\cr c & d\cr}\right)
\right)f\right](x) = |-bx +d|^{-(s+1)}f\left(\frac{ax - c}{-bx + 
d}\right)\, , \]
is now unitary for $0<s<1$. Notice that we now denote by
$\bH_s$ the new Hilbert space
with the inner product $<\cdot \mid A_s(\cdot )>_s$. 

\npar
Define an involution $\t$ on $G$ by
\begin{equation}\label{tforSL2R}
\t\left(\matrix{a & b\cr c & d\cr}\right)
:=\left(\matrix{0 & 1\cr 1 & 0\cr}\right)
\left(\matrix{a & b\cr c & d\cr}\right)
\left(\matrix{0 & 1\cr 1 & 0\cr}\right)
= \left(\matrix{d & c\cr b & a\cr}\right)\, .\end{equation}

The group $H$ is given by
\[ H =\pm \left\{\left. h_t =  \left(\matrix{ \cosh t & \sinh t\cr
\sinh t & \cosh t\cr}\right)\,\right| \,t\in \bbR\right\}\]
and the space $\fq$ is
\[\fq=\left\{\left. q(r,s):= \left(\matrix{ r  & s\cr
-s  & -r\cr}\right)\,\right| \,r,s\in \bbR\right\}\, .\]
Take
\[C := \{q(r,s)\mid r\pm s\ge 0, r\ge 0\} = 
\overline{\conv\{\bbR _+\Ad(H)q(1,0)\}}\, .\]
as a generating cone. The Cartan involution $\th$ is
given by $a\mapsto a^{-t} = waw^{-1}$ and
the corresponding maximal compact subgroup is
$\SO (2)$. 
Define
\[Jf(a) := f(\t (a)w^{-1}) = f (\t (aw))\, .\]
Then $J : \bH_s \pf \bH_s$ intertwines $\p_s$
and $\p_s\circ \t$, and  $J^2=1$. 
In our realization we have
$J(f)(x) = |x|^{-s-1}f(1/x)$
and
\[A_s(J(g))(x) = \int_{-\infty}^\infty 
g(y)|1-xy|^{s-1}\, dy\, .\]
Hence
\[\ip{f}{g}_J = \ip{f}{A_sJg} = \int_{-\infty}^{\infty}\int_{-\infty}^{\infty} 
\overline{f(x)}g(y)|1-xy|^{s-1}\, dy\,dx\, .\]
Let $\bK_0$ be the completion of the space of smooth functions with
compact support in
$I = (-1,1)$. Notice that the above inner product
is defined on $\cC^\infty_c(I)$ for every $s$ as we only integrate
over compact subsets of $(-1,1)$.

\npar
The Bergman kernel for the domain $\{z\in \bbC\mid |z|< 1\}$ is
$h(z,w) = 1 - z\halo{w}$ and it is well known (cf. 
\cite[p. 268]{FK94}) that $h(z,w)^{-\l}$ is a positive definite
kernel function if and only if $\l \ge 0$. As our kernel is
just $h(z,w)^{-(1-s)}$ restricted to the interval $I$ and
and $s < 1$, i.e., $1-s>0$, it follows that $\ip{\cdot }{\cdot }_J$ is positive
definite.

\npar 
We also know (cf. \cite{HO95}) that $S = H\exp C$ is a
closed semigroup and that $\g I \subset I$
and actually $S $ is exactly the semigroup of elements in
$\lSL(2,\bbR)$ that acts by contractions on $I$.
Hence $S$ acts on $\bK$. By a theorem of L\"uscher and Mack \cite{HiNe93,LM75},
the representation of $S$ on $\bK$ extends to a representation
of $G^c$, which in this case is the universal covering of $\SU 
(1,1)$ that is locally isomorphic to $\lSL (2,\bbR)$. We notice that
this defines a representation of $\lSL(2,\bbR)$ if and only
if certain integrality conditions hold; see \cite{JoMo84}.
\npar
We generalize this contstruction to the non-compactly causal
symmetric spaces and in particular to the Cayley-type spaces.
Furthermore we indentify the resulting representation as
an irreducible unitary highest weight representation of the
dual group $G^c$. We restrict ourself to the case of
characters induced from a maximal parabolic subgroup, which
leads to highest weight modules with one-dimensional lowest
$K^c$-type. This is meant as a simplification and not as
a limitation of our method.

\npar Assume now that $G$ is a semidirect product of $H$ and $N$ with
$N$ normal and abelian. Define $\t :G\pf G$ by $\t (hn)=hn^{-1}$. Let
$\p \in \hat{H}$ (the unitary dual) and extend $\p $ to a unitary
representation of $G$ by setting $\p (hn)=\p (h)$. In this case, $G^c$
is locally isomorphic to $G$, and $\p $ gives rise to a unitary
representation $\p ^c$ of $G^c$ by the formula $d\p ^c(X)=d\p (X)$,
$X\in \fh $, and $d\p ^c|_{i\fq} =0$. A special case of this is the
$3$-dimensional Heisenberg group, and the $(ax+b)$-group. In sections
\ref{Diagonal} and \ref{axb}, we show that, if we induce instead a
character of the subgroup $N$ to $G$, then we have $\left( \bK
_0\bftdiagup \bN \right) \tilde{}=\{ 0\} $.

\npar Our approach to the general representation correspondence $\p
\mapsto \p ^c$ is related to the integrability problem for
representations of Lie groups (see \cite{JoMo84}); but the present
positivity viewpoint comes from Osterwalder-Schrader positivity; see
\cite{OsSc73, OsSc75}. In addition the following other papers are
relevant in this connection: \cite{FOS83,Jor86,Jor87,KlLa83,Pra89,Sch86}.

\section{\protect\label{Preliminaries}Preliminaries}
\setcounter{equation}{0}

The setting for the paper is a general Lie group $G$ with
a nontrivial involutive automorphism $\t$.

\begin{Def}\label{ReflectionSymmetric}{\rm A unitary
representation $\p$ acting on a Hilbert space $\bH (\p)$ is 
said to be {\it reflection symmetric} if there is a
unitary operator $J : \bH (\p) \pf \bH(\p)$ such
that
\begin{list}{}{\setlength{\leftmargin}{\customleftmargin} 
\setlength{\itemsep}{0.5\customskipamount}
\setlength{\parsep}{0.5\customskipamount}
\setlength{\topsep}{\customskipamount}
\setlength{\parindent}{0pt}}
\item[\hss\llap{\rm R1)}] ${\displaystyle J^2 = \id}$.
\item[\hss\llap{\rm R2)}] ${\displaystyle J\p (g) = \p (\t (g))J\, ,
\quad g\in G}$.
\end{list}
}
\end{Def}

\npar
If (R1) holds then $\p$ and $\p\circ \t$ are equivalent. 
Furthermore, generally from (R2) we have $J^2\p (g) =
\p (g) J^2$. Thus, if $\p$ is irreducible, then we can always
renormalize $J$ such that (R1) holds.
Let $H = G^\t = \{g\in G\mid \t(g) = g\}$ and let
$\fh$ be the Lie algebra of $H$. Then
$\fh = \{ X\in \fg\mid \t (X) = X\}$. Define
$\fq = \{Y\in \fg\mid \t (Y) = -Y\}$. Then $\fg = \fh \oplus \fq$,
$[\fh,\fq]\subset \fq$ and $[\fq,\fq]\subset \fh$.

\begin{Def}\label{Hyperbolic}{\rm A closed convex cone $C\subset
\fq$ is {\it hyperbolic}
if $C^o\not=\emptyset$ and if $\ad X$ is semisimple with real
eigenvalues for every $X\in C^o$.  }\end{Def}

We will assume the following for $(G, \p,\t, J)$:
\begin{list}{}{\setlength{\leftmargin}{\customleftmargin} 
\setlength{\itemsep}{0.5\customskipamount}
\setlength{\parsep}{0.5\customskipamount}
\setlength{\topsep}{\customskipamount}}
\item[\hss\llap{{\rm PR1)}}] $\p$ is reflection symmetric with reflection $J$.
\item[\hss\llap{{\rm PR2)}}] There is an $H$-invariant hyperbolic cone
$C\subset \fq$ such that $S(C) = H\exp C$ is a closed semigroup and
$S(C)^o = H\exp C^o$ is diffeomorphic to $H\times C^o$.
\item[\hss\llap{{\rm PR3)}}] There is a subspace ${0}\not=
\bK_0\subset \bH(\p)$ invariant under $S(C)$ satisfying the positivity
condition
\[ \ip{v}{v}_J:= \ip{v}{J(v)} \ge 0,\quad \forall v\in \bK_0\, .\]
\end{list}

\begin{Rem}\label{WhatToAssume}{\rm In (PR3) we can always assume
that $\bK_0$ is
closed, as the invariance and the positivity passes over to
the closure. In (PR2) it is only necessary to assume that
$\bK_0$ is invariant under $\exp C$, as one can always replace
$\bK_0$ by $\overline{\left\langle \p (H)\bK_0\right\rangle }$, the
closed space generated by
$\p (H)\bK_0$, which is $S (C)$-invariant, as $C$ is
$H$-invariant. For the exact conditions on the cone
for (PR2) to hold see the orginal paper by J. Lawson
\cite{JL94} or the monograph \cite[pp. 194 ff.]{HiNe93}.
}
\end{Rem}

\npar In some of the examples we will replace (PR2)and
(PR3) by the following
weaker conditions
\begin{list}{}{\setlength{\leftmargin}{\customleftmargin} 
\setlength{\itemsep}{0.5\customskipamount}
\setlength{\parsep}{0.5\customskipamount}
\setlength{\topsep}{\customskipamount}}
\item[\hss\llap{\rm PR2$^{\prime }$)}] $C$ is (merely) some nontrivial
cone in $\fq $.
\item[\hss\llap{\rm PR3$^{\prime}$)}] There is a subspace
$0\not= \bK_0\subset \bH (\p)$ invariant under $H$ and $\exp C$
satisfying the positivity condition from (PR3).
\end{list}
(See Section \ref{Diagonal} for further details.)

\npar
Since the operators $\{\p (h)\mid h\in H\}$ commute with
$J$, they clearly pass to the quotient by
\[\bN := \{v\in \bK_0\mid \ip{v}{Jv} = 0\}\]
and implement unitary operators on $\bK := \left(
\bK_0\bftdiagup \bN\right)\tilde{}$ relative to the inner product
induced by
\begin{equation}\label{E:innerpr}
\ip{u }{ v}_J := \ip{u}{J(v)}\, .
\end{equation}
which will be denoted by the same symbol.
Hence we shall be concerned with passing the operators
$\{\p (\exp Y)\mid Y\in C\}$ to the quotient
$\bK_0\bftdiagup \bN$, and for this we need a basic Lemma.

\npar In general, when $\left( \bK _0,J\right) $ is given, satisfying
the positivity axiom, then the corresponding composite quotient
mapping
\[\bK _0\pf \bK _0\bftdiagup \bN \hookrightarrow \left( \bK _0\bftdiagup
\bN \right) \tilde{}=:\bK \] 
is {\em contractive} relative to the
respective Hilbert norms. The resulting (contractive) mapping will be
denoted $\b $. An operator $\g $ on $\bH $ which leaves $\bK _0$
invariant is said to {\em induce} the operator $\tilde{\g }$ on $\bK $
if $\b \circ \g =\tilde{\g }\circ \b $ holds on $\bK _0$. In general,
an induced operation $\g \mapsto \tilde{\g }$ may not exist; and, if
it does, $\tilde{\g }$ may fail to be bounded, even if $\g $ is
bounded.

\npar This above-mentioned operator-theoretic formulation of
reflection positivity has applications to the Feynman-Kac formula in
mathematical physics, and there is a considerable literature on that
subject, with work by E. Nelson, A. Klein and L.J. Landau, B. Simon,
and W.B. Arveson. Since we shall not use path space measures here, we
will omit those applications, and instead refer the reader to the
survey paper \cite{Arv84} (lecture 4) by W.B. Arveson. In addition to
mathematical physics, our motivation also derives from recent papers
on non-commutative harmonic analysis which explore analytic
continuation of the underlying representations;
see, e.g., \cite{HOO91,Nee94,'O90a,O93,Ol82}.

\section{\protect\label{Basic}A Basic Lemma}
\setcounter{equation}{0}

\begin{Lemma}\label{BasicLemma}
\begin{enumerate}\item[\hss\llap{\rm 1)}] Let $J$ be a period-$2$ unitary
operator on a Hilbert space $\bH$, and let
$\bK_0\subset  \bH$ be a closed subspace such that
$\ip{v}{J(v)}\ge 0$, $v\in \bK_0$. Let
$\g$ be an invertible operator on $\bH$ such that
$J\g = \g^{-1}J$ and which leaves $\bK_0$ invariant and
has $(\g^{-1})^*\g$ bounded on $\bH$. Then $\g$ induces
a bounded operator $\tilde{\g}$ on $\bK = 
\left(\bK_0\bftdiagup \bN\right)\tilde{}$, where
$\bN = \{v\in \bK_0\mid \ip{v}{Jv} = 0\}$, and the
norm of $\tilde{\g}$ relative to the $J$-inner product in $\bK$ 
satisfies
\begin{equation}\label{E:3.1}
\|\tilde{\g}\|\le \|(\g^{-1})^*\g\|^{1/2}_{sp}\, ,
\end{equation}
where $\|\cdot \|_{sp}$ is the spectral radius.
\item[\hss\llap{\rm 2)}] If we have a semigroup $S$ of operators on
$\bH$ satisfying the conditions in {\rm (1)}, then
\begin{equation}\label{E:3.2}
(\g_1\g_2)\tilde{} = \tilde{\g_1}\tilde{\g_2}\, ,
\quad \g_1,\g_2\in S\, .
\end{equation}
\end{enumerate}
\end{Lemma}

{\it Proof\/}: For $v\in \bK_0$, $v\not= 0$, we
have
\begin{eqnarray*}
\|\g (v)\|_J^2 &= & \ip{\g (v)}{J\g (v)}\\
&=& \ip{\g (v)}{\g^{-1}J(v)}\\
&=& \ip{(\g^{-1})^*\g (v) }{J(v)}\\
&=& \ip{(\g^{-1})^*\g(v) }{v}_J\\
&\le& \|(\g^{-1})^*\g (v)\|_J\|v\|_J\\
&\le& \|((\g^{-1})^*\g)^2(v)\|_J^{1/2}\|v\|_J^{1+1/2}\\
&\vdots &\\
&\le & \|((\g^{-1})^*\g)^{2^n}(v)\|_J^{1/2^n}\|v\|_J^{1+1/2+
\cdots + 1/2^n}\\
&\le & \left( \|((\g^{-1})^*\g)^{2^n}\|\|v\|\right) ^{1/2^n}\|v\|_J^{2}\, .
\end{eqnarray*}
Since
$\displaystyle{\lim_{n\to\infty}\|((\g^{-1})^*\g)^{2^n}\|^{1/2^n}
= \|(\g^{-1})^*\g\|_{sp}}$
and
${\displaystyle \lim_{n\to \infty}\|v\|^{1/2^n} = 1}$,
the result follows.

\npar
By this we get
\[\ip{\g (v)}{J\g (v)}  \le \|(\g^{-1})^*\g\|_{sp}\ip{v}{J(v)}\]
which shows that $\g (\bN) \subset \bN$, whence $\g$ passes
to a bounded operator on the quotient $\bK_o\bftdiagup \bN$ and then
also on $\bK$ satisfying the estimate stated in (1).
If both the operators in (\ref{E:3.2}) leave $\bN$ invariant,
so does $\g_1\g_2$ and the operator induced by $\g_1\g_2$ is
$\tilde{\g_1}\tilde{\g_2}$ as stated. \hfill$\Box$\medskip

\begin{Cor}\label{GammaContraction}Let the notation be as above and
assume that $\g$ is unitary on
$\bH$. Then the constant on the right in {\rm (\ref{E:3.1})} is one. 
Hence $\tilde{\g}$ is a contraction on $\bK$.
\end{Cor}

\npar
To understand the assumptions on the space $\bK_0$, i.e.,
positivity and invariance, we include the folowing which
is based on an idea of R.S. Phillips \cite{Phil}.

\begin{Prop}\label{P:3.3} Let $\bH$ be a Hilbert space and let
$J$ be a period-$2$ unitary operator on $\bH$. Let
$S$ be a commutative semigroup of unitary operators
on $\bH$ such that $S = S_+S_-$ with
$S_+ = \{ \g\in S\mid J\g = \g J\}$ and
$S_- = \{\g \in S\mid J\g = \g^{-1}J\}$. Then $\bH$ possesses
a maximal positive and invariant subspace, i.e., a subspace
$\bK_0$ such that
$\ip{v }{J(v)}\ge 0$, $v\in \bK_0$ and $\g \bK_0\subset \bK_0$,
$\g \in S$.
\end{Prop}

{\it Proof\/}: The basic idea is contained in \cite[pp. 386 ff.]{Phil}.
We can represent $\bH$ as $\bL^2(X,m)$ where $X$ is a
Stone space. There is an $m$-a.e.-defined automorphism
$\th : X\pf X$ such that
\[ J(f) = f\circ \th,\quad f\in \bL^2(X,m)\, ,\]
and $S$ is represented by multiplication operators on
$\bL^2(X,m)$. {}From \cite[Lemma 5.1]{Phil}, we know that there are
clopen subsets $A,B$ in $X$ such that if
\[
M_0 = \{x\in X\mid \th(x) = x\}\]
and
\[M_1 = X\bftdiagdown M_0\, ,\]
then $A$ and $B$ are contained in $M_1$,
\begin{eqnarray*}
A\cap B & = & \emptyset\, \\
A\cup B &=& M_1
\end{eqnarray*}
and
\[\th (A) = B\, .\]
Let $\bK_0 := \bL^2(M_0\cup A)$. It is clear that
this is a maximal positive and invariant subspace. The
positivity follows in the following way: If $f$ is supported
in $A$ then $\bar{f}\, f\circ \th = 0$ a.e. Hence for
$f\in \bL^2(M_0\cup A)$,
\begin{eqnarray*}
\ip{f}{J(f)} &=& \int_{M_0} \bar{f}f\circ\th \,dm +
\int_A\bar{f}f\circ \th \,dm\\
&=& \int_{M_0} |f|^2\,dm + \int_A \bar{f}f\circ \th \,dm\\
&=& \int_{M_0}|f|^2\,dm \ge 0\, .
\end{eqnarray*}
This proves the Lemma.\hfill$\Box$\medskip

\begin{Cor}\label{MNoughtZero}If $M_0\subset X$ is of measure 
zero, then the space $\bK$ will be trivial, i.e.,
$\ip{f}{J(f)} = 0$ for all $f\in \bK_0$.
\end{Cor}

\begin{Rem}\label{AbelianSubspace}{\rm Assume that we have (PR1) and (PR2).
Assume
further that we can find an abelian subspace $\fa\subset \fq$ such
that $C^o = \Ad (H)(C^o\cap \fa)$. Let $S_A = \exp (C^o\cap
\fa)$. Then $S_A$ is an abelian semigroup, so one can use Proposition
\ref{P:3.3} to construct a maximal positive and invariant subspace for
$S_A$. But in general we can't expect this space to be invariant under
$S$.  }
\end{Rem}

\npar
We read off from the basic Lemma that

\begin{Prop}\label{Contractive}Let $\p$ be a unitary representation
of $G$. Assume that $(\t, J, C, \bK_0)$
satisfies the conditions {\rm (PR1)}, {\rm (PR2${}^\prime$)} and
{\rm (PR3${}^\prime$)}.
If $Y\in C$ then $\p (\exp Y)$ induces a contractive selfadjoint
operator $\tilde{\p}(\exp Y)$ on $\bK$.
\end{Prop}

{\it Proof\/}: If $Y\in C$ then $\p (\exp Y)\bK_0\subset \bK_0$
and
$\p (\exp Y)$ is unitary on $\bH(\p)$. Thus
\begin{eqnarray*}
\ip{\p (\exp Y)u}{J(v)} & = & \ip{u}{\p (\exp (-Y))J(v)}\\
&=& \ip{u}{J(\p (\exp Y)v)}\, ,
\end{eqnarray*}
proving that $\p (\exp Y)$ is selfadjoint in
the $J$-inner product. Since $\p (\exp Y)$ is unitary
on $\bH (\p )$
\[\| \p(\exp Y ) \| = \|\p (\exp Y)\|_{sp} = 1\, ,\]
and the contractivity property follows. \hfill$\Box$\medskip

\begin{Cor}\label{HolomorphicRPlus}Let $\p$ be a unitary representation
of $G$ such that $(\t, J, C, \bK_0)$
satisfies the conditions 
{\rm (PR1)}, {\rm (PR2${}^\prime$)} and
{\rm (PR3${}^\prime$)}.
Then for
$Y\in C$ there is a selfadjoint operator 
$d\tilde{\p}(Y)$ in $\bK = \left(\bK_0\bftdiagup \bN\right)\tilde{}$ with
spectrum contained in $(-\infty , 0]$ such that
\[\tilde\p (\exp (tY)) = e^{t d\tilde{\p}(Y)},\quad t\in \bbR _+\]
is a contractive semigroup on $\bK$.
Furthermore the following hold:
\begin{enumerate}
\item[\hss\llap{\rm 1)}] $t\mapsto e^{td\tilde{\p}(Y)}$ extends to a
continuous  map
$z \mapsto e^{zd\tilde{\p}(Y)}$ on $\{z \in \bbC\mid
\Re (z)\ge 0\}$ holomorphic on the
open right half-plane and
such that $e^{(z + w)d\tilde{\p}(Y)} =e^{zd\tilde{\p}(Y)}e^{wd\tilde{\p}(Y)}$.
\item[\hss\llap{\rm 2)}] If $Y\in C^o$ then the above map is holomorphic in an
open neighborhood of  $\{z \in \bbC\mid
\Re (z)\ge 0\}$.
\item[\hss\llap{\rm 3)}] There
exists a one-parameter group of
unitary operators
\[\tilde{\p}\left(\exp (itY)\right):= e^{itd\tilde{\p }(Y)},
\quad t\in \bbR\]
on $\bK$.
\end{enumerate}
\end{Cor}

{\it Proof\/}: The last statement follows by the spectral
theorem. By construction
$\{\tilde{\p }(\exp (tY))\mid t\in \bbR _+\}$ is a semigroup
of selfadjoint contractive operators on $\bK$. The
existence of the operators $d\tilde{\p }(Y)$ as stated
then follows from a general result in operator theory; see,
e.g.,
\cite{Fr80} or \cite{KlLa81}.\hfill$\Box$\medskip

\begin{Cor}\label{Id}Let the situation be as in the last
corollary. If $Y\in C\cap -C$ then $e^{td\tilde{\p}(Y)}
= \id$ for all $t\in \bbR _+$. In particular
$d\tilde{\p }(Y) = 0$ for every $Y\in C\cap -C$.
\end{Cor}

{\it Proof\/}: This follows as the spectrum
of $d\tilde{\p}(Y)$ and $d\tilde{\p }(-Y)$ is contained
in $(-\infty , 0]$.\hfill$\Box$\medskip

\section{\protect\label{LM}The L\"uscher-Mack Theorem}
\setcounter{equation}{0}

We use reference \cite{HiNe93} for the L\"uscher-Mack Theorem,
but \cite{FOS83}, \cite{GoJo83}, \cite{Jor86}, \cite{Jor87},
\cite{JoMo84}, \cite{KlLa83}, \cite{LM75}, and \cite{Sch86} 
should also be mentioned in this connection.

\npar
Let $\p$, $C$, $\bH(\p)$, $J$ and $\bK_0$ be as before.
We have proved that the operators
\[\{\p (h \exp (Y))\mid h\in H, Y\in C\}\]
pass to the space $\bK = \left(\bK_0\bftdiagup \bN\right)\tilde{}$ such that
$\tilde{\p}(h)$ is unitary on $\bK$, and
$\tilde{\p}(\exp Y)$ is contractive and selfadjoint on $\bK$.
As a result we arrive at selfadjoint operators $d\tilde{\p}(Y)$
with spectrum in $(-\infty, 0]$ such that for $Y\in C$,
$\tilde{\p }(\exp Y) = e^{d\tilde{\p}(Y)}$ on $\bK$.
As a consequence of that we notice that 
\[ t\mapsto e^{t d\tilde{\p}(Y)}\]
extends to a continuous map on $\{z\in \bbC\mid \Re (z) \ge 0\}$
holomorphic on the open right half plane $\{z\in\bbC\mid \Re (z) >0\}$.
Furthermore,
\[e^{(z+w)d\tilde{\p }(Y)} = e^{zd\tilde{\p }(Y)}e^{wd\tilde{\p }(Y)}\, .\]
As $\bK$ is a unitary $H$-module we know that the $H$-analytic vectors
$\bK^\omega(H)$ are dense in $\bK$. Thus $\bK_{oo} := S
(C^o)\bK^\go(H)$ is dense in $\bK$. We notice that for $u\in \bK_{oo}$
and $X\in C^o$ the function $t\mapsto \tilde{\p }(\exp tX)u$ extends
to a holomorphic function on an open neighborhod of the right
half-plane. This and the Campbell-Hausdorff formula are among the main
tools used in proving the following Theorem of L\"uscher and Mack
\cite{LM75}.  We refer to \cite[p. 292]{HiNe93} for the proof. Our
present use of Lie theory, cones, and semigroups will follow standard
conventions (see, e.g.,
\cite{FHO93,Hel62,JL94,WaI72,Yos91}): the
exponential mapping from the Lie algebra $\fg $ to $G$ is denoted
$\exp $, the adjoint representation of $\fg $, $\ad $, and that of $G$
is denoted $\Ad $. If $\p $ is a representation of $G$, its
differential is denoted $d\p $, e.g., $d(\Ad )=\ad $. Recall that if
$\p $ is infinite-dimensional, then $d\p $ is a representation by
unbounded operators on $\bH (\p )$, but the analytic vectors and the
$C^\infty $-vectors form dense domains for $d\p $; see \cite{Nel59,Pou92}.

\begin{Satz}[L\"uscher-Mack]\label{LuscherMack}
Let $\r $ be a strongly continuous contractive representation
of $S (C)$ on the Hilbert space $\bH $ such that
$\r (s)^* = \r (\t (s)^{-1})$. Let
$G^c$ be the connected, simply connected  Lie group with Lie
algebra $\fg^c = \fh\oplus i\fq$. Then there exists a
continuous unitary representation $\r^c : G^c
\pf {\rm U}(\bH )$ extending $\r $ such that for the
differentiated representations
$d\r$ and $d\r^c$ we have:
\begin{enumerate}
\item[\hss\llap{\rm 1)}] $d\r^c(X) = d\r (X)\, \quad \forall X\in \fh$.
\item[\hss\llap{\rm 2)}] $d\r^c(iY) = i\, d\r (Y)\, \quad \forall Y\in C$.
\end{enumerate}
\end{Satz}

\npar
We apply this to our situation to get:
\begin{Satz}\label{PiCIrreducible}Assume that $(\p , C, \bH,J)$ satisfies
{\rm (PR1)--(PR3).} Then the following hold:

\begin{enumerate}
\item[\hss\llap{\rm 1)}] $S(C)$ acts via $s\mapsto \tilde{\p}(s)$ by
contractions on $\bK$.
\item[\hss\llap{\rm 2)}] Let $G^c$ be the simply connected Lie group with Lie
algebra $\fg^c$. Then there exists a unitary representation
$\tilde{\p}^c$ of $G^c$ such that 
$d\tilde{\p}^c(X) = d\tilde{\p}(X)$ for $X\in \fh$ and
$i\, d\tilde{\p}^c(Y) = d\tilde{\p}(iY)$ for $Y\in C$.
\item[\hss\llap{\rm 3)}] The representation $\tilde{\p}^c$ is irreducible
if and only if $\tilde{\p}$ is irreducible.
\end{enumerate}
\end{Satz}

{\it Proof\/}: (1) and (2) follow by the L\"uscher-Mack theorem
and Proposition 3.6, as the resulting representation of
$S$ is obviously continuous.

\npar
(3) Let $\bL$ be a $G^c$-invariant subspace in $\bK$.
Then $\bL$ is $\tilde{\p}(H)$ invariant. Let $Y\in C^o$,
$u\in \bL^\go$ and $v\in \bL^\perp$. Define
$f : \{z\in \bbC\mid \Re (z) \ge 0\} \pf \bbC$ by
\[ f(z) := \ip{v}{ e^{z d\tilde{\p}(Y)}u}_J\, .\]
Then $f$ is holomorphic in $\{z\in \bbC\mid \Re (z) >0\}$,
and $f(it) = 0$ for every (real) $t$. Thus $f$ is identically
zero. In particular $f(t) = 0$ for every $t>0$. Thus
\[0 = \ip{v}{e^{td\tilde{\p}(Y)}u}_J = \ip{v}{\tilde{\p}(\exp tY)u}_J\, .\] 
As $S^o = H\exp C^o$
it follows that $\tilde{\p}(S^o)(\bL^\go )\subset (\bL^\perp)^\perp = \bL$.
By continuity we get $\tilde{\p}(S)\bL \subset \bL$. Thus
$\bK$ is reducible as an $S$-module.

\npar
The other direction follows in exactly the same way.\hfill$\Box$\medskip

\npar
Let $(\p , C,\bH,J)$ be as in the last theorem. 
To identify the resulting representation 
$\tilde{\p}^c$ of $G^c$ some facts about holomorhic
representations of semigroups and highest weight
representations are needed. We refer to
\cite[Chap. 7]{HO95} and the references therein, in
particular \cite{Nee94}, for further references. Define
\[W(\tilde{\p }^c) := \{ X\in \fg^c\mid \forall u\in \bK^\infty :
i\ip{u}{\tilde{\p }^c(X)u}_J\le 0\}\] where $\bK^\infty$ denotes the
$\cC^\infty$-vectors for $G^c$.  Then $W(\tilde{\p }^c)$ is a closed
$G^c$-invariant cone in $\fg^c$.  $W (\tilde{\p }^c)$ is non-trivial as
$-iC\subset W (\tilde{\p }^c)$.  Thus $W(\tilde{\p }^c)$ will always
contain the $-\t$-stable and $G$-invariant cone generated by $-iC$,
i.e.  $-i\Ad (G)C$, but in general $W(\tilde{\p }^c)$ is neither
generating nor pointed. It even does not have to be
$-\t$-invariant. In fact, the Lie algebra of the $(ax+b)$-group, and
the Heisenberg group, does not have {\em any} pointed, generating,
invariant cone.

\begin{Lemma}\label{KernelPiC}$W(\tilde{\p }^c)\cap -W(\tilde{\p
}^c) = \ker (\tilde{\p }^c)$.
\end{Lemma}

{\it Proof\/}: This is obvious from the spectral
theorem.\hfill$\Box$\medskip

\begin{Lemma}\label{IdealGC}$\fg_1^c := W(\tilde{\p }^c)-W(\tilde{\p }^c)$
is an ideal in
$\fg^c$. Furthermore, $[\fq,\fq]\oplus i\fq\subset \fg_1^c$.
\end{Lemma}

{\it Proof\/}: Let $X\in \fg^c$. Then, as $W (\tilde{\p }^c)$ is
invariant by construction, we conclude that 
\[e^{t\ad (X)}\left(W
(\tilde{\p }^c)-W(\tilde{\p }^c)\right) \subset W(\tilde{\p
}^c)-W(\tilde{\p }^c),\, t\in \bbR\, .\]
By differentiation at $t=0$, it
follows that $[X,\fg_1^c]\subset \fg_1^c$. This shows that $\fg_1^c$
is an ideal in $\fg^c$. The last part follows as $C$ is generating (in
$\fq$).\hfill$\Box$\medskip

\begin{Rem}\label{TauStable}{\rm It is not clear if $\fg_1^c$ is
$\t$-stable. To get a $\t$-stable subalgebra one can replace $W
(\tilde{\p }^c)$ by the cone generated by $-\Ad (G) C\subset W
(\tilde{\p }^c)$ or by the maximal $G$- and $-\t$-stable cone
$W(\tilde{\p }^c)\cap (-\t (W(\tilde{\p }^c)))$ in $W(\p^c)$.  }
\end{Rem}

\npar
Let $W$ be a $G_1^c$ invariant cone in $\fg_1^c$. We define
$\cA (W)$ to be the set of equivalence classes of
irreducible unitary representations $\r$ of $G_1^c$ with
$W (\r )\subset W$.
 
\begin{Satz}\label{HighestWeight}Assume that the analytic subgroup
$G^c_1$ of $G^c$ corresponding
to $\fg_1^c$ is closed in $G$ and that $W(\tilde{\p }^c)$ is pointed.
Then $\tilde{\p }^c|_{G_1^c}$ is a direct integral of irreducible
representations in $\cA (W)$.
\end{Satz}

{\it Proof\/}: As $G_1^c$ is closed in $G$ it follows that
$\tilde{\p }^c|_{G_1}$ is a continuous unitary representation of
$G_1^c$. Furthermore $W (\tilde{\p }^c|_{G_1}) = W(\tilde{\p }^c)$. The theorem
follows now from the theorem of Neeb and Olshanskii \cite{Nee94}, that
an injective representation $\r$ with
$W (\r )$ pointed and generating is a direct integral of
representations from $\cA (W (\r ))$ (cf. \cite{Nee94}).\hfill$\Box$\medskip

\section{\protect\label{SSS}Examples of semisimple symmetric spaces}
\setcounter{equation}{0}

We will now generalize the example from the Introduction
to a class of semisimple Lie groups. For that we
recall some facts about non-compactly
causal or ordered semisimple symmetric spaces. We include some
ideas of the proofs to make the text more
self contained. For more information we refer to  \cite{HO95,'O90b}.
An additional source of inspiration for the present chapter is the
following series of papers: \cite{Nel59,'OO89a,OO96}.

\npar
Let $G\bftdiagup H$ be a semisimple symmetric space and
let $\t$ be the corresponding involution. {\it We
will assume that $G\bftdiagup H$ is irreducible}. Let $\th$ be
a Cartan involution on $G$ commuting with $\t$. Then
\begin{eqnarray*}
\fg &=& \fh \oplus \fq\\
&=& \fk\oplus \fp\\
&=& \fh_k\oplus \fh_p \oplus \fq_k\oplus \fq_p
\end{eqnarray*}
where a subscript denotes the intersection with the
corresponding subspace of $\fg$.
Let $L$ be a Lie group and $\bV$ an $L$-module. We denote by
$\bV^L$ the subspace of $L$-fixed points in $\bV$.

\begin{Def}\label{NonCompactlyCausal}
{\rm The irreducible symmetric space $G\bftdiagup H$ is called {\it
non-com\-pact\-ly causal} (NCC) if $\fq_p^{H\cap K}\not= \{0\}$.
}\end{Def}

\begin{Rem}{\rm A NCC-space is also a $K_\e$-space in the
sense of \cite{OsS1,OsS2}. 
}
\end{Rem}
 
\npar
If $G\bftdiagup H$ is NCC then $\fq_p^{H\cap K}$ is one-dimensional
and there exists an element $X^0\in \fq_p^{H\cap K}$ such that
$\fh_k\oplus \fq_p = \fz_{\fg}(X^0)$. We can
normalize $X^0$ such that
$\ad X^0$ has eigenvalues $0$, $1$, and $-1$.
Let $\fa := \bbR X^0$, $\fn = \{X\in \fg\mid
[X^0,X] = X\}$, and
$\halo{\fn} = \{X\in \fg\mid [X^0,X] = -X\} = \th (\fn)
= \t (\fn)$.
We also define
\[\fm := \{X\in \fz_\fg(X^0)\mid B(X,X^0) = 0\}\]
where $B$ is the Killing form of $\fg$. Then
\[\fp_{\rm max} := \fm \oplus \fa\oplus \fn\]
is a maximal parabolic subalgebra of $\fg$.

\npar
Assume from now on
that $G\subset G_{\bbC}$ where $G_{\bbC}$ is the simply
connected, connected Lie group with  Lie algebra
$\fg_{\bbC}$. We will also assume that
$H = G^\t$. Then $H\cap K = Z_K(X^0)$.
Let $A := \exp \fa$, $N := \exp \fn$ and
$\halo{N} := \exp \halo{\fn}$.
Let $M_o$ be the analytic subgroup of $G$ corresponding
to $\fm$ and let $M = (H\cap K)M_o$. Then
$M$ is a closed and $\t$-stable subgroup of $G$, $M\cap A
= \{1\}$
and $MA = Z_G(A)$.
Let $P_{\rm max} := N_G(\fp_{\rm max})$.
Then $P_{\rm max} = MAN$.
We have
$\fg = \fh + \fp_{\rm max}$. The differential
of the map $(h,p)\mapsto hp$ is given by
$\Ad (p)(X + \Ad (p^{-1})Y)_{hp}$, $X\in \fh$
and $Y\in \fp_{\rm max}$, and this is surjective, as
$\Ad (p^{-1})\fp_{\rm max} = \fp_{\rm max}$.

\begin{Lemma}\label{HPminOpenG}$HP_{\rm min}$ is open in $G$ and
contained in $\halo{N}P_{\rm max}$.
\end{Lemma}

{\it Proof\/}: That $HP_{\rm max}$ is open in $G$ follows by
the above discussion (for the general case see
\cite{Ma82}).
The proof of the second statement can be found in
\cite{HO95,'O90b}. The idea is to use a maximal
set of strongly orthogonal roots to 
reduce this to $\lSL (2,\bbR)$-calculations as
we will explain in a moment.
\hbox to 0.1cm{}\hfill$\Box$\medskip

\npar Let $\fa_q$ be a maximal abelian subalgebra of $\fp$ containing
$X^0$. Then $\fa_q\subset \fq_p$ and $\fa_q$ is maximal abelian in
$\fq$. Let $\D$ be the set of roots of $\fa_q$ in $\fg$. Then $\D =
\D_0\cup \D_+\cup \D_-$, where $\D_0 = \{\a\in \D\mid \a (X^0) = 0\}$,
$\D_\pm = \{\a\in \D\mid \a (X^0) = \pm\, 1\}$.  Choose a positive
system $\D_0^+$ in $\D_0$ and let $\D^+ = \D_0^+\cup \D_+$. Two roots
$\a,\b$, $\a\not= \pm \b$ are called strongly orthogonal if $\a \pm
\b$ is not a root. Choose a maximal set of strongly orthogonal roots
$\g_1 < \g_2 < \cdots < \g_r$ in $\D_+$ such that $\g_r$ is the
maximal root in $\D_+$, $\g_{r-1}$ is the maximal root in $\D_+$
strongly orthogonal to $\g_r$, $\g_{r-2}$ is the maximal root in
$\D_+$ strongly orthogonal to $\g_r$ and $\g_{r-1}$, etc. Choose
$H_j\in \fa_q$ such that $\left\langle \g_i,H_j\right\rangle =
2\d_{ij}$ and $H_j\in [\fg_{\g_j},\fg_{-\g_j}]$.  Choose $X_j\in
\fg_{\g_j}$ such that with $X_{-j} := \t (X_j) = - \th (X_j)$ we have
$H_j = [X_j,X_{-j}]$. In the case of $\lsl (2,\bbR)$ the involution is
given by \[\t\left(\matrix{a & b\cr c & d\cr}\right) :=\left(\matrix{0
& 1\cr 1 & 0\cr}\right) \left(\matrix{a & b\cr c & d\cr}\right)
\left(\matrix{0 & 1\cr 1 & 0\cr}\right) = \left(\matrix{d & c\cr b &
a\cr}\right)\] as in the Introduction. In this case
\[H_1 = \left(\matrix{1 & 0 \cr 0 & -1\cr}\right),
\quad X_1=\left(\matrix{0 & 1 \cr 0 & 0\cr}\right)
\quad\mbox{\rm and}\quad
X_{-1} = \left(\matrix{0 & 0 \cr 1 & 0\cr}\right)\, .\]
Define a homomorphism $\f_j : \lsl (2,\bR)\pf \fg$ by
\begin{eqnarray*}
\left(\matrix{1 & 0 \cr 0 & -1\cr}\right) & \pf & H_j\\
\left(\matrix{0 & 1 \cr 0 & 0\cr}\right) & \pf & X_j\\
\left(\matrix{0 & 0 \cr 1 & 0\cr}\right)
&\pf & X_{-j}
\end{eqnarray*}
As the roots $\g_j$ are strongly orthogonal we get
$[\Im (\f_j),\Im (\f_i)] = \{0\}$ if $i\not= j$. As $\lSL (2,\bbC)$
is simply connected the  homomorphisms $\f_j$ integrate
to homomorphisms $\lSL (2,\bbC) \pf G_{\bbC}$, also denoted
by $\f_j$, such that $\f_j (\lSL (2,\bbR)) \subset G$ and
such that $\f_j$ intertwines the Cartan involution and
the above involution $\t$ on $\lSL (2,\bbR)$ with
the corresponding involutions on $G$.

\npar
The following lemma follows from the maximality of
the set of strongly orthogonal roots; see also 
\cite[Lemma 2.3]{'OO88a}:

\begin{Lemma}\label{MaximalAbelian}Let $\fa_h = \bigoplus_{j=1}^r
\bbR (X_j+ X_{-j})$.
Then $\fa_h\subset \fh_p$, and $\fa_h$ is maximal abelian
in $\fh_p$.
\end{Lemma}

%

\npar
Let $\log := (\exp|_{\halo{\fn}})^{-1} : \halo{N}\pf \halo{\fn}$.
Define $\z : \halo{N}P_{\rm max}\bftdiagup P_{\rm max} \pf
\halo{\fn}$ by
\[\z (\halo{n}P_{\rm max}) = \log (\halo{n})\, .\]
We notice that $\z (hx) = \Ad (h)\z (x)$ for $h\in H\cap K$.
We also notice that $H\cap P_{\rm max} = H\cap K$.

\begin{Lemma}\label{TanhSum}Let $h = \exp \sum_{j=1}^r t_j(X_j+X_{-j})\in 
A_h := \exp \fa_h$. Then
\[hP_{\rm max} = (\exp \sum_{j=1}^r \tanh t_j X_{-j})P_{\rm max}\, .\]
In particular $\z (hP_{\rm max}) = \sum_{j=1}^r \tanh t_j X_{-j}$.
\end{Lemma}

{\it Proof\/}: Assume first that $G = \lSL(2,\bbR)$.  Then $h = h_t =
\left(\matrix{\cosh t & \sinh t\cr \sinh t & \cosh t\cr}\right)$.  By
(\ref{NbarPmax}), we have $h_t \in \halo{N}P_{\rm max}$, and $\z
(h_tP_{\rm max} ) = \tanh t X_{-1}$.  Let $t_1,\ldots ,t_r\in \bbR$:
then

\newlength{\expsumetc} \settowidth{\expsumetc}{$\displaystyle \exp
\sum_{j=1}^r t_j(X_j + X_{-j})P_{\rm max}$}
\addtolength{\expsumetc}{-2em}
\begin{eqnarray*}
&\ &\mbox{\makebox[\expsumetc][r]{\hss $\displaystyle \exp
\sum_{j=1}^r t_j(X_j + X_{-j})P_{\rm max}$}}\\ 
&\ & \qquad \qquad =
\f_1(h_{t_1})\cdots \f_r(h_{t_r})P_{\rm max}\\ 
&\ & \qquad \qquad =
\f_1\left(\left(\matrix{ 1 & 0 \cr \tanh t_1 & 1 \cr}\right)
\right)\cdots \f_r\left(\left(\matrix{ 1 & 0 \cr \tanh t_r & 1
\cr}\right) \right)P_{\rm max}\\ 
&\ & \qquad \qquad =(\exp
\sum_{j=1}^r \tanh (t_j)X_{-j}) P_{\rm max}\, .
\end{eqnarray*}
{F}rom this the lemma now follows. \hfill$\Box$\medskip 

\begin{Satz}\label{OmegaConvex}Let $\gO = \Ad
(H\cap K)\{\sum_{j=1}^rt_jX_{-j}\mid
\forall j : -1<t_j<1\}$. Then $\gO$ is convex,
\[ HP_{\rm max} = (\exp \gO)P_{\rm max}\, ,\]
and $\z$ induces an isomorphism $H\bftdiagup H\cap K \simeq \gO$.
\end{Satz}

{\it Proof\/}: The convexity will follow from
Lemma \ref{T:compdom}.
That $(\exp \gO)P_{\rm max} \subset HP_{\rm max}$
follows from the fact that
$\exp t_jX_{-j}P_{\rm max}\subset
HP_{\rm max}$ by $\lSL (2,\bbR)$-reduction.
Let $h\in H$: then $h$ can be written as
$h = k_1ak_2$, with $k_1,k_2\in H\cap K$ and
$a = \exp \sum t_j(X_j+X_{-j})\in A_h$. As
$H\cap K \subset P_{\rm max}$, it follows that
\begin{eqnarray*}
hP_{\rm max}& =& k_1aP_{\rm max}\\
&=& k_1\exp \sum_{j=1}^r \tanh (t_j)X_{-j} P_{\rm max}\\
&=&\exp ( \Ad (k_1) \sum_{j=1}^r \tanh (t_j)X_{-j}) P_{\rm max}\\
&\in & \exp \gO P_{\rm max} \subset \halo{N}P_{\rm max}
\end{eqnarray*}
Thus $HP_{\rm max}\subset \exp (\gO)P_{\rm max}$.\hfill$\Box$\medskip 

\npar
The maximal compactly embedded subalgebra $\fk^c$ in
$\fg^c$ corresponding to the Cartan involution
$\th^c = \th\t$ has center $i\fa$ and $\fz_{\fg^c}(i\fa) = \fk^c$.
It follows that $G^c\bftdiagup K^c$ is a bounded symmetric domain
and that $\t$ induces an anti-holomorphic involution
on $G^c\bftdiagup K^c$, i.e.~a
conjugation. The real form of $G^c\bftdiagup K^c$ corresponding
to this conjugation is exactly $H\bftdiagup H\cap K$, (see
\cite{HJ75,HJ78} for classification). In the classical
notation of Harish-Chandra (cf. \cite{He84}) we
have $\fp^- = \halo{\fn}_{\bbC}$. Thus $G^c\bftdiagup K^c$ can
be realized as a bounded symmetric domain $\gO_{\bbC}$
in $\halo{\fn}_{\bbC}$. Let $\s $ be the conjugation
of $\fg_{\bbC}$ with respect to $\fg$. Then
$\s|\fg^c = \t|\fg^c$. Thus the conjugation given by
$\t$ on $\gO_{\bbC}$ is also realized by $\s$. We now have:

\begin{Lemma}\label{T:compdom} Let $\gO_{\bbC}$ be the
bounded convex circular realization of $G^c\bftdiagup K^c$
in $\halo{\fn}_{\bbC}$. Then
\[\gO = \gO_{\bbC}^\s = \{X\in \gO_{\bbC}\mid
\s (X) = X\}\, .\]
\end{Lemma}

\npar
We also notice the following for later use:

\begin{Lemma}\label{ConjugationCoincides}Denote the conjugation
of  $\fg_{\bbC}$ with
respect to $\fg^c$ by $\s^c$. Then $\s^c$ coincides with
the conjugate linear extension $\t\circ \s$
of $\t$ to $\fg_{\bbC}$.
\end{Lemma}

{\it Proof\/}: We have
$\{X\in \fg_{\bbC}\mid \t\s (X) = X\} = \fh\oplus i\fq = \fg^c$.
Hence the lemma.\hfill$\Box$\medskip

\npar
Let
\[S(H,P_{\rm max}) : = \{g\in G\mid g H \subset HP_{\rm max}\}
\, .\]
Then $S(H,P_{\rm max})$ is a closed semigroup invariant
under $s\mapsto s^\sharp := \t(s)^{-1}$.
For $g\in G$ and $X\in \halo{\fn}$ such that
$g\exp X\in\halo{N}P_{\rm max}$ define
$g\cdot X\in \halo{\fn}$ and $a(g,X)\in A$ 
by
\[g\exp X \in \exp (g\cdot X)Ma(g,X)N\, .\]
In particular we have $g\cdot X = \z (g\exp X)$, where $\z
:H\bftdiagup H\cap K\simeq \gO $ is introduced in the Theorem
\ref{OmegaConvex} (see also Lemmas
\ref{MaximalAbelian}--\ref{TanhSum}), and it follows that the elements
$a(g,X)$ are defined for $g\in S(H,P_{\rm max}) $ and $X\in \gO$.  The
map $(g,X)\mapsto g\cdot X$ transfers the canonical action on
$G\bftdiagup P_{\rm max}$ restricted to the open set $HP_{\rm
max}\bftdiagup P_{\rm max}$ to $\gO$. We have

\begin{Lemma}\label{L:5.8}\begin{enumerate}
\item[\hss\llap{\rm 1)}] Let $s,r\in S(H,P_{\rm max})$ and
$X\in \gO$. Then
$(sr)\cdot X = (s\cdot (r\cdot X))$ and
$a(sr,X) = a(s,r\cdot X)a(r,X)$.
\item[\hss\llap{\rm 2)}] Let $g = ma \in MA$ and $X\in \halo{\fn}$.
Then $g\exp X \in \halo{N}P_{\rm max}$,
$g \cdot X = \Ad (g)X$, and $a(g,X) = a$.
\item[\hss\llap{\rm 3)}] Let $C$ be an  $H$-invariant pointed
and generating cone in
$\fq $ containing $X^0$. Then $S=H\exp C$ is a closed
semigroup acting on $\gO $ by contractions. Furthermore $H\times
C^o\ni (h,X)\mapsto h\exp X \in S^o$ is a diffeomorphism.
\item[\hss\llap{\rm 4)}] $S(H,P_{\rm max})\subset HP_{\rm max}$.
\end{enumerate}
\end{Lemma}

{\it Proof\/}: Let $s,r$ and $X$ be as in the lemma. Then
on the one hand
\[(sr)\exp X = \exp ((sr)\cdot X) m (sr,X)a(sr,X)n(sr,X)\]
for some $m(sr,X)\in M$ and $n(sr,X)\in N$. On the other hand, using
the notation
\begin{eqnarray*}
n &=& n(r,X)\, ,\\
n_1 &=& n(s,r\cdot X)\, ,\\
n_2 &=& [(m(r,X)a(r,X))^{-1}n(s,r\cdot X)m(r,X)a(r,X)]n(r,X)
\end{eqnarray*}
($n,n_1,n_2\in N$), we have
\begin{eqnarray*}
(sr)\exp (X) &=& s(r\exp (X))\\
&=& s \exp (r\cdot X) m (r,X)a (r,X)n\\
&=& \exp (s\cdot (r\cdot X)) m(s,r\cdot X)a(s,r\cdot X)n_1
m(r,X)a(r,X)n_1\\
&=& \exp (s\cdot (r\cdot X))m(s,r\cdot X)m (r,X)
a(s,r\cdot X)a(r,X)n_2
\, .
\end{eqnarray*}
This proves (1).

\npar
(2) This follows from 
\[g\exp X = \exp (\Ad (g)X)g\]
and the fact that $MA$ normalizes $\halo{N}$.

\npar (3) Let $\fp $ and $\fq $ be as described before Definition
\ref{NonCompactlyCausal} above, and let $C$ be a pointed and
generating $H$-invariant cone in $\fq$ such that $C^o\cap \fp \not=
\emptyset$.  Then by \cite{HO95,'O90b}
\begin{equation}\label{c0=}
 C^o = \Ad (H)(C^o\cap \fa_q)\, .
\end{equation}

Let $X\in C^o\cap \fa$ and $Y\in \gO$. Then $Y = \sum Y_{-\a}$, with
$\a \in \D_+$ and $Y_{-\alpha }\in \fg _{-\alpha }$. Therefore
\[\exp (X)\cdot Y = \Ad (\exp X) Y
= \sum _{\alpha \in \D _+}e^{-\a (X)}Y_{-\a}\]
As $\a (X)> 0$ we see that $\exp (X)\cdot Y \in \gO$. This also shows that
$\exp (C^o\cap \fa)$ acts by contractions on $\gO$. 
Let now $s\in S^o$. Then $s = h\exp X = h \exp (\Ad (h_1)X_1)$
with $h,h_1\in H$, $X\in C^o$ and $X_1\in C^o\cap \fa$.
Let $Y\in \gO$. Then
\[s\cdot Y = hh_1\cdot \left( \exp X_1\cdot \left( \left( h_1^{-1}
\right) \cdot Y\right) \right) \in \gO\, .\]
It follows that $S$ acts by contractions on $\gO$.

\npar
(4) Apply $S(H,P_{\rm max})$ to $eP_{\rm max}$, $e $ the
identity in $G$.\hfill$\Box$\medskip

\begin{Lemma}\label{Simple}Let $t>0$ and $Y\in \gO$. Then
$\exp tX^0\in S$ and $\exp tX^0\cdot Y = e^{-t}Y$.
\end{Lemma}

\npar
We also notice the following sharpening of (3) in
Lemma \ref{L:5.8} (cf. \cite{HiNe93} and \cite{HO95}):

\begin{Lemma}\label{Sharpen} Let $C =C_{\rm max}$
be the maximal poirefear
nted generating cone in $\fq$ containing
$X^0$. Then the following holds:
\begin{enumerate}
\item[\hss\llap{\rm 1)}]  $C^o\cap \fa =\left\{ X\in \fa _q\mid \forall \alpha
\in \D _+:\alpha (X)>0\right\} $.
\item[\hss\llap{\rm 2)}] 
$S(H,P_{\rm max} )= H\exp C_{\rm max}$.
\end{enumerate}
\end{Lemma}

\npar
We need to fix the normalization of measures before
we discuss the representations that we will use.
Let the measure $da$ on $A$ be given by
\[\int_A f(a)\,da = \frac{1}{\sqrt{2\p}}\int_{-\infty}^\infty\,
f(a_t)\,dt, \quad a_t = \exp 2tX_0\, .\]
We fix the Lebesgue measure $dX$ on $\halo{\fn}$ such that for
$d\halo{n} = \exp (dX)$ we have
\[\int_{\halo{\fn}}a(\halo{n})^{-2\r}\,d\halo{n} = 1\, .\]
Here $\r (X) = \frac{1}{2}\Tr (\ad (X))|_{\fn}$ as usual,
and $a(g)$, $g\in G$, is determined by $g\in KMa(g)N$.
The Haar measure on compact groups will always be normalized
to have total measure one. The measure on $N$ is $\th (d\halo{n})$.
Let us fix a Haar measure $dh$ on $H$. 
Then we can normalize the invariant measure on
$G$ such that for $f\in \cC_c(G)$, $\Supp (f)\subset HP_{\rm max}$,
we have
\[\int_G f(g)\,dg = \int_H\int_A\int_N f(han)a^{2\r}\,dn\,da\,dh\, .\]
The invariant measure $d\dot{x}$ on $G\bftdiagup H$ is then given by
\[\int_G f(x)\,dx = \int_{G\bftdiagup H}\int_H f(xh)\,dh\,d\dot{x},\quad
f\in \cC_c(G)\]
and similarly for $K\bftdiagup H\cap K$.
We fix the Haar measure on $M$ such that
$dg = a^{2\r}\,dk\,dm\,da\,dn$. 

\begin{Lemma}\label{L:Int} Let the measures be normalized as above. Then
the following hold:
\begin{enumerate}
\item[\hss\llap{\rm 1)}] Let $f\in \cC_c(\halo{N}MAN)$. Then
\[\int_Gf(g)\,dg = \int_{\halo{N}}\int_M\int_A\int_N
f(\halo{n}man)a^{2\r}\,d\halo{n}\,dm\,da\,dn\, .\]
\item[\hss\llap{\rm 2)}] Let $f\in \cC_c(\halo{N})$. For $y\in \halo{N}
MAN$ write
$y = \halo{n}(y)m_{\halo{N}}(y)a_{\halo{N}}(y)n_{\halo{N}}(y)$.
Let $x\in G$. Then
\[\int_{\halo{N}}f(\halo{n}(x\halo n))a_{\halo{N}}(x\halo n
)^{-2\r}\,d\halo n =
\int_{\halo{N}} f(\halo n)\,d\halo{n}\, .\]
\item[\hss\llap{\rm 3)}] Write, for $g\in G$, $g=k(g)m(g)a(g)n(g)$
according to $G = KMAN$.
Let $h\in \cC (K\bftdiagup H\cap K)$.
Then
\[\int_{K\bftdiagup H\cap K} h(\dot k) \,d\dot k =
\int_{\halo N} h(k(\halo n) H\cap K) a(\halo n)^{-2\r}\,d\halo n \, .\]
\item[\hss\llap{\rm 4)}] Let $h\in \cC (K\bftdiagup H\cap K)$ and
let $x\in G$.
Then
\[\int_{K\bftdiagup H\cap K}f (k(xk)H\cap K)a(xk)^{-2\r}\,d\dot{k} =
\int_{K\bftdiagup H\cap K}
f(\dot k)\,d\dot k\]
\item[\hss\llap{\rm 5)}] Assume that $\Supp (f)\subset H\bftdiagup
H\cap K\subset K\bftdiagup H\cap K$. Then
\[\int_{K\bftdiagup H\cap K} f(\dot k)\,d\dot k = \int_{H\bftdiagup
H\cap K}f (k (h)H\cap K)
a(h)^{-2\r}\,d\dot h\, .\]
\item[\hss\llap{\rm 6)}] Let $f\in \cC_c(\halo N)$. Then
\[\int_{\halo N} f(\halo n)\,d\halo n = 
\int_{H\bftdiagup H\cap K} f(\halo{n}(h))a_{\halo N}(h)^{-2\r}\,d\halo n\, .\]
\item[\hss\llap{\rm 7)}] For $x \in HP_{\rm max}$ write
$x = h(x)m_H(x)a_H(x)n_H(x)$ with $h(x)\in H$,
$m_H(x)\in M$, $a_H(x)\in A$ and $n_H(x)\in N$. Let
$f\in \cC_c^\infty (H\bftdiagup H\cap K)$ and let
$x\in G$ be such that $xH P_{\rm max} \subset
HP_{\rm max}$. Then
 \[\int_{H\bftdiagup H\cap K}f (h(xh)H\cap K)a_H(xh)^{-2\r}\,d\dot{h} =
\int_{H \bftdiagup H\cap K}
f(\dot h)\,d\dot h\]

\end{enumerate}
\end{Lemma}

\npar
{\it Proof\/}: Up to normalizing constants this can be
found in \cite{'O87}. Let us show that the constant in (1) is
equal to 1. Choose $c>0$ such that
\[c\int_Gf(g)\,dg = \int_{\halo{N}}\int_M\int_A\int_N
f(\halo{n}man)a^{2\r}\,d\halo{n}\,dm\,da\,dn\, .\]
Let $\f \in \cC_c(P_{\rm max})$ such that
$\int_{MAN}\f(man)a^{2\r}\,dm\,da\,dn =1$
and $\f (mp) = \f (p)$ for every $m\in M\cap K = H\cap K$ 
and every $p\in P_{\rm max}$. Define
$f\in \cC(G)$ by $f (kman) = \f (man)$.
Then
\begin{eqnarray*}
c &=& c\int_K\int_M\int_A\int_Nf(kman)a^{2\r}\,dk\,dm\,da\,dn\\
&=& c\int_G f(g)\,dg\\
&=& \int_{\halo{N}}\int_M\int_A\int_Nf(\halo n man)a^{2\r}\,d\halo n\,dm
\,da\,dn\\
&=& \int_{\halo{N}}\int_M\int_A\int_Nf(k(\halo n)m(\halo 
{n}) a(\halo{n})n(\halo{n}) man)a^{2\r}\,d\halo n\,dm
\,da\,dn\\
&=&\int_{\halo{N}}\int_M\int_A\int_Nf(k(\halo n)
man)a(\halo{n})^{-2\r}a^{2\r}\,d\halo n\,dm
\,da\,dn\\
&=& \int_{\halo N}a(\halo n)^{-2\r}\,d\halo n = 1
\end{eqnarray*}
This proves (1). The other claims are proved in a similar
way. \hfill$\Box$\medskip

\npar
Let us now go over to the representations that we are going to
use. We identify $\fa_{\bbC}^*$ with $\bbC$ by
\[\fa^*_{\bbC}\ni\n \mapsto 2\n (X^0)\in \bbC\, .\]
Then $\r$ corresponds to $\dim \fn$.
For $\n \in \fa^*_{\bbC}$, let $\cC^\infty(\n)$ be the
space of $\cC^\infty$-functions $f : G\pf \bC$ such
that, for $a_t = \exp t (2X^0)$,
\[f(gma_tn) = e^{-(\n +\r)t}f(g) = a_t^{-(\n + \r)}f(g)\, .\]
Define an inner product on $\cC^\infty(\n)$ by
\[\ip{f}{g}_\n := \int_K\, \overline{f(k)}g(k)\,dk =
\int_{K\bftdiagup H\cap K}\, \overline{f(k)}g(k)\,d\dot k  \, .\]
Then $\cC^\infty(\n)$ becomes a pre-Hilbert space.
We denote by $\bH(\n)$ the completion of $\cC^\infty(\n )$.
Define $\p (\n)$ by
\[[\p (\n )(x)f](g) := f(x^{-1}g),\quad x,g\in G,\quad f
\in \cC^\infty(\n)\, .\]
Then $\p (\n ) (x)$ is bounded, so it extends to a
bounded operator on $\bH(\n)$, which we denote by
the same symbol and
$\p (\n)$ is a continuous representation of
$G$ which is unitary if and only if $\n \in i\bbR$. 
Furthermore $\bH(\n)^\infty = \cC^\infty(\n)$ (cf. \cite{Pou92}).
We can realize $\bH(\n)$ as $\bL^2(K\bftdiagup H\cap K)$
and as $\bL^2(\halo{N}, a (\halo{n})^{2\Re (\n)}\,d\halo{n})$
by restriction (see Lemma \ref{L:5.14}).
In the first realization the representation $\p (\n)$ becomes
\[[\p (\n)(x)f](k) = a(x^{-1}k)^{-\n-\r}f(k(x^{-1}k))\]
and in the second
\[[\p (\n)(x)f](\halo n) =
a_{\halo N}(x^{-1}\halo{n})^{-\n-\r}f(\halo{n}(x^{-1}\halo n))\, .\]
The following is well known, but for completeness we include
the proof:

\begin{Lemma}\label{L:invpar} The pairing
\begin{eqnarray*}
\bH(\n)\times \bH(-\halo{\n})
\ni (f,g)\mapsto \ip{f}{g} _\n &:= &\int_K\overline{f(k)}g(k)\,dk\\
&=& \int_{K\bftdiagup H\cap K} \overline{f(k)}g(k)\,d\dot{k}
\end{eqnarray*}
is $G$-invariant, i.e.
\[\ip{\p (\n)(x)f}{g}_\n = \ip{f}{\p(-\halo{\n})(x^{-1})g}_\n\, .\]
\end{Lemma}

{\it Proof\/}: Let $x\in G$ and $k\in K$. Then
$x(x^{-1}k) = k$, which implies that
\begin{eqnarray*}
k &= & xk(x^{-1}k)a(x^{-1}k)n(x^{-1}k)\\
&= &k(xk(x^{-1}k))a(xk(x^{-1}k))n(xk(x^{-1}k))a(x^{-1}k)n(x^{-1}k)\\
&=&k(xk(x^{-1}k))a(xk(x^{-1}k))a(x^{-1}k)n
\, .
\end{eqnarray*}
for some $n = \in N$.
Thus $k(xk(x^{-1}k)) = k$
and $a(xk(x^{-1}k)) = a(x^{-1}k)^{-1}$. Using those relations, and
Lemma \ref{L:Int},
we get:
\begin{eqnarray*}
\ip{\p (\n)(x)f}{g}_\n &=& \int_{K\bftdiagup H\cap
K}\overline{f(x^{-1}k)}g(k)\,d\dot k\\
&=& \int_{K\bftdiagup H\cap
K}\overline{f(k(x^{-1}k))}a(x^{-1}k)^{-(\halo{\n}+\r)}
g(k(xk(x^{-1}k)))\,d\dot k\\
&=&\int_{K\bftdiagup H\cap K}\overline{f(k)}\left[a(xk)^{-(-\halo{\n}+\r)}
g(k(xk))\right]\,d\dot k\\
&=&\ip{f}{\p (-\n)(x^{-1})g}_\n\, .
\end{eqnarray*}
This proves the lemma.\hfill$\Box$\medskip

\begin{Rem}\label{ImaginaryNu}{\rm We notice that if $\n $ is purely imaginary,
i.e., $-\halo{\n} = \n$, the above shows that
$(\p (\n),\bH (\n))$ is then unitary.
}
\end{Rem}

\begin{Lemma}\label{L:5.14}{\ }
\begin{enumerate}
\item[\hss\llap{\rm 1)}] The restriction map induces an isometry of
$\bH (\n)$ onto $\bL^2(\halo N, a(\halo n)^{2\Re \n}\,d\halo n )$.
\item[\hss\llap{\rm 2)}] On $\halo{N}$ the invariant pairing
$\ip{\cdot }{\cdot }_\n$ is given by
\[ \ip{f}{g}_\n = \int_{\halo N}\overline{f(\halo n)}g(\halo n)
\,d\halo n\, ,\quad f\in \bH (\n), \, g\in \bH(-\halo \n) \, .\]
\item[\hss\llap{\rm 3)}] Let $\bH_H(\n)$ be the closure of $\{f\in
\cC^\infty(\n)\mid \Supp (f)\subset HP_{\rm max}\}$.  Then
$\bH_H(\n)\ni f\mapsto f|_H\in \bL^2(H\bftdiagup H\cap
K,a(h)^{2\r}\,d\dot h)$ is an isometry.
\item[\hss\llap{\rm 4)}] Let $f\in \bH (\n), \, g\in \bH(-\halo \n)$ 
and assume that $\Supp (fg)\subset HP_{\rm max}$. Then
\[\ip{f}{g}_{\n} = \int_{H\bftdiagup H\cap K}\overline{f(h)}g(h)\,d\dot h\, .\]
\end{enumerate}
\end{Lemma}

{\it Proof\/}: (1)
We have
$k(\halo{n}) = \halo{n}a(\halo n)^{-1}n$, $\halo{n}\in 
\halo{N}$. By Lemma \ref{L:Int} we get

\[\int_{K\bftdiagup H\cap K}\overline{f(k)}g (k)\,d\dot k
= \int_{\halo N} \overline{f(\halo n)}g(\halo n)
a (\halo n)^{\halo{\n} + \m}\,d\halo n\, ,
\quad f\in \bH (\n), \, g\in \bH(\m)\, .\]

\npar
(2)--(4) follow immediately from Lemma \ref{L:Int}.\hfill$\Box$\medskip

\npar
Let us assume, from now on, that there exists an element
$w\in N_{K}(\fa)$ such that $\Ad (w) (X^0) = -X^0$ on $\fa$.
Let us remark the following for later use:

\begin{Lemma}\label{L:tw} Let $w\in K$ be such that $\Ad (w)|_\fa = -\id$.
Then
$w^2\in H\cap K$ and there is a $m\in H\cap K$ such
that $\t (w) = w^{-1}m$.
\end{Lemma}

{\it Proof\/}: As $\Ad (w^2)|_\fa = \id$ we get
$w^2 \in M\cap K = H\cap K$. Let $X\in \fa$.
Then 
\[X = \t (\Ad (w)X) = \Ad (\t (w))(\t (X)) = - \Ad (\t(w))X\, .\]
Hence $\Ad (\t (w)w)X = X$. Thus
$\t (w)w =: m \in M\cap K = H\cap K$.
It follows that $\t (w) = w^{-1}(wmw^{-1})$. The claim follows
as $wMw^{-1} = M$.\hfill$\Box$\medskip

\npar
We recall that $G\bftdiagup H$ is of {\it Cayley type} if
$\fh$ has a one-dimensional center contained in
$\fh_p$. This is the case if and only if $G\bftdiagup K$ is a tube-type
domain $G\bftdiagup K\simeq \bbR^n+i\gO$, where $\gO$ is an
open self-dual cone isomorphic to $H\bftdiagup H\cap K$.
Thus $G\bftdiagup H$ is locally isomorphic to one of the 
following spaces (where we denote by the subscript $+$ the group
of elements having positive determinant):
$\Sp (n,\bbR )\bftdiagup \GL (n,\bbR)_+$, $\SU (n,n)\bftdiagup \GL (n,\bbC)_+$,
$\SO^*(4n)\bftdiagup \SU^*(2n)\bbR_+$, $\SO(2,k)\bftdiagup \SO (1,k-1)\bbR_+$  
and $E_{7(-25)}\bftdiagup E_{6(-26)}\bbR_+$.

\begin{Lemma}\label{CayleyType1}Assume that $G\bftdiagup H$ is of Cayley type.
Let
\[ w= \f_1\left(\left( \matrix{ 0 & 1\cr -1 & 0 \cr}\right)\right)\cdots
\f_r\left(\left( \matrix{ 0 & 1\cr -1 & 0 \cr}\right)\right)
= \exp \left(\frac{\p}{2}\sum_{j=1}^r X_j + \th (X_j)\right)\, .\]
Then $\Ad (w)|_\fa = -\id$.
\end{Lemma}

{\it Proof\/}: As $G\bftdiagup H$ is of Cayley type,
$X^0 = \frac{1}{2}\sum_{j=1}^rH_j$. The claim follows now
by simple $\lsl (2,\bbR)$-calculation.\hfill$\Box$\medskip

\npar
We also recall the following lemma from \cite{HO95,'O90b}:

\begin{Lemma}\label{CayleyType2}Assume that $G\bftdiagup H$ is
of Cayley type. Let
$Y^0\in \fh_p$ be such that
$\fz_{\fg}(Y^0) = \fh$ and such that
$\spec (\ad Y^0) = \{0,1,-1\}$. 
Then 
$\bc := \Ad (\exp \frac{\p i}{2}Y^0)$ defines
a Lie algebra isomorphism $\fg\pf \fg^c$
such that
\begin{enumerate}
\item[\hss\llap{\rm 1)}] $\bc|_\fh = \id_{\fh}$.
\item[\hss\llap{\rm 2)}]
Let
$\fq^+ := \{X\in \fq\mid [Y^0,X] = X\}$
then $\bc|_{\fq^+} = i\id$.
\item[\hss\llap{\rm 3)}] Let
$\fq^- := \{X\in \fq\mid [Y^0,X] = -X\}$. Then
$\bc|_{\fq^-} = -i \id$.
\item[\hss\llap{\rm 4)}] $\fq = \fq^+\oplus \fq^-$.
\end{enumerate}
\end{Lemma}

{\it Proof\/}: 
That $\bc : \fg\pf \fg^c$ is an
isomorphism follows from (1)--(4).
(1)--(3) follow directly. For (4) notice
that $\ad Y^0$ maps $\fq$ into $\fq$.
As the centralizer of $Y^0$ is exactly $\fh$ it
follows that $\ad Y^0 :\fq \pf \fq$ is an
isomorphism and that $\fq$ is the direct
sum of the eigenspaces of $\ad Y^0$ for the
eigenvalues $1$ and $-1$. {}From that
the claim follows.\hfill$\Box$\medskip

\npar
Assume now that $\fh$ is one of the
Lie algebras $\sp (n,\bbR )$, $\su (n,n)$,
$\so^*(4n)$, $\so (2,k)$  
and $\fe_{7(-25)}$. Let $\fg = \fh_{\bbC}$
and let $G_{\bbC}$ be the simply connected connected
group with Lie algebra $\fg$. Let $\t :\fg\pf \fg$
be the conjugation with respect to $\fh$. Denote the
corresponding real analytic involution $G\pf G$ by the
same letter. Then it is well known that
$G^\t = H$ is connected. We refer to \cite[Example 1.2.2]{HO95}.

\begin{Lemma}\label{LocallyIsomorphic}Assume that $\t : \fg \pf \fg$
is the conjugation
with respect to the real form $\fh$. Then 
\begin{enumerate}
\item[\hss\llap{\rm 1)}] $\fg^c\simeq \fh\times \fh$ and
$G^c$ is locally isomorphic to $H\times H$.
\item[\hss\llap{\rm 2)}] Under this isomorphism the involution $\t$ corresponds
to $\t (X,Y) = (Y,X)$, i.e. $\fh$ corresponds to the
diagonal in $\fg^c$.
\item[\hss\llap{\rm 3)}] Let $\tilde{H}$ be the simply connected
Lie group with Lie algebra $\fh$.
Then $G^c$ is $\tilde{H}\times \tilde{H}$ and
$\t$ is given by $\t (a,b) = (b,a)$. In particular
$(G^c)^\t = \{ (a,a)\mid a\in \tilde{H}\} \simeq \tilde{H}$
and
$G^c\bftdiagup \tilde{H}\ni (a,b)\tilde{H}\mapsto ab^{-1}\in \tilde{H}$
is an isomorphism.
\end{enumerate}
\end{Lemma}

\npar
Notice that in this case we can construct, using the
strongly orthogonal roots, commuting homomorphisms
$\f_j^{\bbC} : \lSL (2,\bbC) \pf G$ such that
actually $\f_j^{\bbC}(\SU (1,1)) \subset H$
and $X^0 = \frac{1}{2}\sum_j\f_j^{\bbC}\left(\left(
\matrix{1 & 0\cr
0 & -1\cr}\right)\right)$. Using this
homomorphism instead of $\f_j$ we get: 
\begin{Lemma}
Let 
\[ w= \prod _j\f_j^{\bbC}\left(\left( \matrix{ 0 & 1\cr -1 & 0
\cr}\right)\right)\, .\]
Then $\Ad (w)(X^0) = - X^0$.
\end{Lemma}

{\it Proof\/}: This follows again by simple $\lsl (2,\bbR)$-calculation
as $X^0 = \frac{1}{2}\sum_{j=1}^r H_j$.\hfill$\Box$\medskip

\npar
For $\Re(\n)$ ``big'' we can construct an intertwining
operator $A(\n) : \bH (\n)\pf \bH (-\n)$ (cf.
\cite{KS80,Wal92}) by
\begin{equation}
[A(\n)f](x) := \int_{\halo{N}}f(xw\halo{n})\,d\halo{n} 
\end{equation} 
Let us show that $A (\n)f \in \bH (-\n)$. For that
let $x\in G$, $man\in MAN$. Then
\begin{eqnarray*}
[A(\n )f](xman) & =&  \int_{\halo{N}}f(xmanw\halo{n})\,d\halo{n} \\
&=&\int_{\halo{N}}f(xw (w^{-1}mw)
a^{-1}(w^{-1}nw)\halo{n})\,d\halo{n}\\
&=& a^{\n+\r}\int_{\halo{N}}f(xw(a^{-1}\halo{n}a))\,d\halo{n}\\
&=& a^{-(-\n+\r)}\int_{\halo{N}}f(xw\halo{n})\,d\halo{n}   
\end{eqnarray*}
Here the third equation follows by the facts that
$w^{-1}Nw = \halo N$, $w^{-1}Mw = M$, and $M$ acts unimodularly on
$\halo N$. The last equation follows by
\[\int_{\halo N} f(a^{-1}\halo n a)\,d\halo n =
a^{-2\r }\int_{\halo N} f(\halo n )\,d\halo n\, .\]
The intertwining property is obvious. 

\npar
The map $\n \mapsto A(\n)$ has an analytic continuation
to a meromorphic function on $\fa^*_{\bbC}$. Because of
Lemma \ref{L:invpar} we can define a new invariant bilinear form
on $\cC^\infty_c(\n)$ by
\[ \ip{f}{g} := \ip{f}{A(\n) g}_\n\, .\]
If there exists a (maximal) constant $R>0$ such that
the invariant bilinear form $\ip{\cdot }{\cdot }$ is
positive definite for $|\n| < R$, we call the resulting
unitary representations {\it the complementary series}.
Otherwise we set $R=0$. We have the following
results from \cite{OZ95}:

\begin{Lemma}\label{CayleyConstant}For the Cayley-type symmetric spaces the
constant $R$ is given by
\begin{eqnarray*}
\SU(n,n)\, :\,  R &=& \left\{\matrix{n &,& n\, \mbox{odd}\cr
0  &,& n\, \mbox{even}\cr}\right.\\
\SO^*(4n)\, : \, R &=& n\cr
\Sp (n,\bbR) \, :\,  R &=& \left\{\matrix{n/2 &,& n\, \mbox{even}\cr
0  &,& n\, \mbox{odd}\cr}\right.\\
\SO_o(n,2) \, :\,  R &=& \left\{\matrix{0 &,& n\equiv 0 \quad
\mbox{\rm mod}\quad  4\cr
1  &,& n\equiv 1,3  \quad\mbox{\rm mod}\quad 4\cr
2 &,& n\equiv 2 \quad \mbox{\rm mod}\quad 4\cr}\right.\\
E_{7(-25)}\, :\, R & = & 3
\end{eqnarray*}
\end{Lemma}

\npar In the cases where $\ip{\cdot}{A(\n )\,\cdot }_\n $ is positive
definite we complete $\cC_c^\infty (\n )$ with respect to this new
inner product, but denote the resulting space by the same symbol
$\bH(\n )$ as before.

\begin{Lemma}\label{PhiUnimodular}$w^{-1}\t (\halo{N})w =\halo{N}$ and
$\f : \halo{N}\ni \halo{n}\mapsto w^{-1}\t (\halo{n})w\in \halo{N}$
is unimodular.
\end{Lemma}

{\it Proof\/}: The first claim follows as $\Ad (w)$ and $\t$ acts by $-1$
on $\fa$ and thus maps $N$ onto $\halo{N}$ and
$\halo{N}$ onto $N$. The second follows as
we can realize 
$\f^2$ by conjugation by an element in $M\cap K$.\hfill$\Box$

\begin{Lemma}\label{UnitaryIsomorphism}For $f\in \bH(\n)$ let
$J(f)(x) := f(\t (xw))$. Then the
following properties hold:
\begin{enumerate}
\item[\hss\llap{\rm 1)}] $J(f)(x) = f(\t(x)w^{-1})$.
\item[\hss\llap{\rm 2)}] $J(f)\in \bH(\n)$ and $A(\n )J=JA(\n )$.
\item[\hss\llap{\rm 3)}] $J :\bH(\n)\pf \bH(\n)$ is
an unitary isomorphism.
\item[\hss\llap{\rm 4)}] $J^2 = \id$.
\item[\hss\llap{\rm 5)}] For $x\in G$ we have
$J\circ \p(\n)(x) = \p (\n)(\t (x))\circ J$.
\end{enumerate}
\end{Lemma}

{\it Proof\/}: (1) This follows from
Lemma \ref{L:tw}, as $f$ is $M$-right invariant.

\npar
(2) Let $x\in G$ and $man\in P_{\rm max}$. 
By (1) we get:
\begin{eqnarray*}
J(f)(xman) &=& f(\t(x)\t(m)a^{-1}\t (n)w^{-1})\\
&=& f(\t (x)w^{-1}(w\t (m)w^{-1})a (w\t (n)w^{-1}))\\
&=& a^{-(\n + \r)}f(\t (x)w^{-1})\, ,
\end{eqnarray*}
as $\t (M) = M$, $w^{-1}Mw = M$, and
$w^{-1}Nw = \t (N) = \halo{N}$.  For $\Re (\n )$ ``big'' we have
\begin{eqnarray*}
A(\n )[Jf](x) &=& \int Jf(xw\halo{n})\,d\halo{n}\\
&=& \int f(\t (xw\halo{n})w^{-1})\,d\halo{n}\\
&=& \int f(\t (x)\t (w)\t (\halo{n})w^{-1})\,d\halo{n}
\end{eqnarray*}
{F}rom Lemma \ref{L:tw} it follows easily that $\t (w)=m_1w$, for some
$m_1\in M$. Thus by Lemma
\ref{PhiUnimodular}:

\begin{eqnarray*}
A(\n )[Jf](x) &=& \int f(\t (x)m_1w\t (\halo{n})w^{-1})\,d\halo{n}\\
&=& \int f(\t (x)m_1\halo{n})\,d\halo{n}\\
&=& \int f(\t (x)w^{-1}w\halo{n})\,d\halo{n}\\
&=& J[A(\n )f](x)\, .
\end{eqnarray*}
The claim now follows by analytic continuation.

\npar
(3) Using that $\t (dk) = dk$ and that $K$ is unimodular it
follows by direct calculation and (2) that
$J^* = J$. That $J$ is a unitary isomorphism follows now
by (4).

\npar
(4) This follows as $\t^2 = \id $ and $w^2\in H\cap K$.

\npar
(5) Let $x,y\in G$. Then
\begin{eqnarray*}
J[\p (\n)(x)f](y) &=& [\p (\n)(x)f](\t (yw))\\
&=& f(x^{-1}\t (yw))\\
&=& f(\t (\t (x)^{-1}yw))\\
&=& [\p (\n)(\t (x))(Jf)](y)\, ,
\end{eqnarray*}
which is exactly what we wanted to prove.\hfill$\Box$\medskip

\npar
Notice that, even if the individual operators $A(\n)$ and
$J$ do not exists, it is always possible to define
the composite operator $A(\n)J$ by
\[ [A(\n)J](f)(x) := \int_{\halo{N}}\, f(\t(x)\halo{n})\, d\halo{n}\]
for $\Re \l$ ``big'' and then by analytic continuation for
other parameters. By simple calculation we get:

\begin{Lemma}\label{L:AnuJ} Assume that $G\bftdiagup H$ is
non-compactly causal. Then $A (\n)J $ intertwines
$\p (\n)$ and $\p (-\n)\circ \t$ if $A(\n)J$ has no pole at $\n$.
\end{Lemma}

\npar
The next theorem shows that the intertwining operator $A(\n)J$ is
a convolution operator with kernel 
$y, x\mapsto a_{\halo{N}}(\t (y)^{-1}x)^{\n -\r}$. The reflection
positivity then reduces to the problem to determine those
$\n$ for which this kernel is positive semidefinite.

\begin{Satz}\label{ANuJf}Let $f
\in \cC^\infty (\n)$. Then
\[ [A(\n)J](f)(\halo{n})
= \int_{\halo N} f(x)a_{\halo N}(\t (\halo{n})^{-1}x)^{\n - \r}\,dx\, .\]
If $\Supp (f) \subset HP_{\rm max}$ then for $h\in H$
\[[A(\n )J](f)(h) = \int_{H\bftdiagup H\cap K}f (x)a_{\halo N}
(h^{-1}x)^{\n - \r}\,d\dot x\, .\]
\end{Satz}

{\it Proof\/}: We may assume that $\n$ is big enough such
that the integral defining $A(\n)$ converges. The
general statement follows then by analytic continuation.
We have
\begin{eqnarray*}
[ A(\n)J]f (\halo n) &=& \int_{\halo{N}} Jf (\halo{n}wx)\,dx\\
&=& \int f (\t (\halo n) w^{-1}\t (x)w)\,dx\\
&=&\int f (\t (\halo n) x) \,dx\\
&=& \int f(\halo{n}(\t (\halo n) x))
a_{\halo N}(\t (\halo n) x)^{-(\n + \r)}\,dx\, .
\end{eqnarray*}
Now $a_{\halo N}(\t (\halo n) x) = a_{\halo N}(\t (\halo{n})^{-1}
\halo{n}(\t (\halo n) x))^{-1}$. By Lemma \ref{L:Int}
we get 
\begin{eqnarray*}
[ A(\n )J] f(\halo n) & =& \int f(\halo{n}(\t (\halo n) x))
a_{\halo N}(\t (\halo{n})^{-1}\halo{n}(\t (\halo n) x))^{(\n - \r)}
a_{\halo N}(\t (\halo n) x)^{- 2 \r} \,dx\\
& =& \int f( x)
a_{\halo N}(\t (\halo{n})^{-1}x)^{(\n - \r)} \,dx
\end{eqnarray*}
The second statement follows in the same way.\hfill$\Box$\medskip
\begin{Cor}\label{fgJ}Let $f,g\in \bC^\infty(\n)$. Then
\[\ip{f}{g}_J = \int_{\halo{N}}\int_{\halo{N}}
\overline{f(x)}g(y)a_{\halo N}(\t (x)^{-1}y)^{\n - \r}\,dx\,dy\, .\]
If $f$ and $g$ both have support in $HP_{\rm max}$, then
\[\ip{f}{g}_J = \int_{H\bftdiagup H\cap K}\int_{H\bftdiagup H\cap K}
\overline{f(h)}g(k)a_{\halo N}(h^{-1}k)^{\n - \r}\,dh\,dk\, .\]
\end{Cor}

\npar In Theorem \ref{PR} bellow, we use this for describing
the representations for which the corresponding $J$ sesquilineqar
form $\ip{\cdot}{\cdot}_J$ is positive semidefinite on the
space of functions supported on $HP_{\rm max}$.

\npar
Assume that $G\bftdiagup H$ is non-compactly causal.
Let $\cC^\infty_c(\gO)$ be the space of $\cC^\infty$-functions
on $\halo{N}$ with compact support in
$\gO$. 
We view this as the subspace in  
$\cC^\infty (\n)$ consisting of functions
$f$, such that $\Supp (f)\subset HP_{\rm max}$
and $\Supp (f|\halo N)$ is compact.
Then $\ip{f}{g}_J$ is defined for every $f,g\in \cC_c^\infty (\gO)$.
In particular we can form the form $\ip{\cdot}{\cdot}_J$
in all cases.

\begin{Lemma}\label{CSInvariant}Suppose that
$G\bftdiagup H$ is non-compactly causal. Let $s\in S$
and $f\in \cC_c^\infty (\gO)$.
Then $\p (\n)(s)f \in \cC_c^\infty (\gO)$, i.e., $\cC_c^\infty (\gO)$
is $S$-invariant.
\end{Lemma}

{\it Proof\/}: 
Let $f\in \cC^\infty_c(\gO)$ and $s\in S$.
Then $\p (\n)(s)f(x) = f(s^{-1}x) \not= 0$ only if
$s^{-1}x\in \Supp (f) \subset HP_{\rm max}$. Thus
$\Supp (\p (\n)(s)f) \subset s\Supp (f) \subset s HP_{\rm max}
\subset HP_{\rm max}$.\hfill$\Box$\medskip

\npar Let $(\p ,\bH )$ be an admissible representation of $G^c$ and
let $\bH _{K^c}$ be the space of $K^c$-finite elements in $\bH $. For $\d
\in \hat{K^c}$ let $\bH (\d )$ be the subspace of $K^c$-finite vectors of
type $\d $, i.e.,
\[\bH (\d )=\bigcup _{T\in \Hom _{K^c}(\bH _\d ,\bH )}T(\bH _\d )\, ,\]
where $\bH_\d$ is the representation space of $\d$.

\begin{Def}\label{HWRWLKCTX}{\rm $(\p ,\bH )$ is called a {\em
highest-weight representation} of $G^c$ (with respect to $\D^+$) if there
exists a $\d \in \hat{K^c}$ such that
\begin{enumerate}
\item[\hss\llap{\rm 1)}] $d\p (\fn_\bbC) \bH (\d )=0$,
\item[\hss\llap{\rm 2)}] $d\p (U(\halo{\fn} ))\bH (\d )=\bH_{K^c}$.
\end{enumerate}
Notice that the multiplicity of $\d $ in $\p $ is one if $\p $
is irreducible. We call $\d$ for the minimal $K^c$-type of
$\p$.
}\end{Def}

\npar Assume that
$G\bftdiagup H$ non-compactly causal.
By the theorem of Moore (cf. \cite{He78}) we know that the roots
in $\D_+ $ restricted to the span of $H_{1}, \ldots ,H_{r}$,
are given by $\pm\frac{1}{2}(\g_i + \g_j)$, $1\le i\le j\le r$
and possibly $\frac{1}{2}\g_j$. The
root spaces for $\g_j$ are all one-dimensional and the root spaces
$\fg_{\pm \frac{1}{2}(\g_i+\g_j)}$, $1\le i<j\le r$, have all the
common dimension $d$.

\begin{Satz}[Vergne-Rossi,~Wallach]\label{S:VRW} Assume
that $G\bftdiagup H$ is non-compactly causal and
that $G^c$ is simple. Let $\l_0\in \fa^*$ be such
that $<\l_0, H_{r}> = 1$. Let $\g = <\l_0,X^0>$ and let

\[L_{\rm pos}:= -\frac{\g (r-1)d}{2}\, .\]

Then the following
holds:
\begin{enumerate}
\item[\hss\llap{\rm 1)}] For
$\n - \r < L_{\rm pos}$ there exists a irreducible unitary highest weight
representation $(\r_\n, \bK_\n)$ of $G^c$ with
one-dimensional minimal $K^c$-type $\n -\r$.  
\item[\hss\llap{\rm 2)}] If $G\bftdiagup H$ is of Cayley-type, then
$\g = r$. Furthermore 
$\n \le L_{\rm pos}$ if and only if $\n \le  r$.
\end{enumerate}
\end{Satz}
 
{\it Proof\/}: (1) By \cite[pp. 41--42]{VR76} (see also \cite{NW79})
$(\r_\n,\bK_\l)$ exists if $<\n - \r, H_{r}> \le - \frac{(r-1)d}{2}$.
But $\n - \r = <\n -\r,H_r>\l_0$. Hence
$<\n - \r, 2X^0> = <\n - \r,H_r><\l_0,2X^0> = \g <\n -\r, H_r>$,

\npar
(2) If $\GH$ is of Cayley type then
$2 X^0 = \sum_{j=1}^r H_j$ and
$\g_j = \g_r - \sum n_\a\a$, $\a \in \D^+_0$, $n_\a\ge 0$.
Thus $<\n -\r, X^0> = r<\n - \r, H_r>$. We also have
(cf. \cite{OO96b})
\[\r = \frac{1}{2}\left( 1 + \frac{(r-1)d}{2}\right)(\g_1 + \cdots
+ \g_r)\, .\]
{}From this the theorem follows.
\hfill$\Box$\medskip

\npar
Let us state this more
explicitely for
the Cayley-type spaces to compare
the existence of $(\r_\n,\bK_\n)$ to
the existence of the complementary
series, cf. Lemma \ref{CayleyConstant}:

\begin{Lemma}\label{CayleyConstantLps}For the Cayley-type symmetric
spaces the highest weight representation $(\r_\n,\bK_\n)$ exists
for $\n$ in the following half-line:

\begin{eqnarray*}
\SU(n,n)\, &:&\,  \n \le n.\\
\SO^*(4n)\, &:& \, \n \le 2n\\
\Sp (n,\bbR) \, &:&\,  \n \le n\\
\SO_o(n,2) \, &:&\,  \n \le 2\\
E_{7(-25)}\, &:&\, \n \le 3\, .
\end{eqnarray*}
In particular we have that
$(\r_\n,\bK_\n)$ is defined for
$\n \in [-R,R]$.
\end{Lemma}

\begin{Rem}{\rm Let us remind the reader that we have only described
here the continuous part of the unitary spectrum. There are also
finitely many discrete points, the socalled {\it Wallach set}, giving
rise to unitary highest weight representations.
}
\end{Rem}

\npar
Let us still assume that $G^c$ is simple.
Let $(\r_\n,\bK_\n)$ be as above. Let
$u\in \bK_\l (\l - \r)$, $\| u\| =1$.
Let $H^c = (G^c)^\t$. Then $H^c$ is connected
\cite{'O87}. Let $\tilde{H}$ be the universal
covering of $H^c $ and $H_o$. We notice that
\[H^c\bftdiagup H^c \cap K^c = H \bftdiagup H\cap K
= H_o\bftdiagup \tilde{H}/\tilde{H}^\th\, .\]
Denote the restriction of $\r_\n$ to $H^c$ by
$\r_{\n,H}$. We can then lift $\r_{\n,H}$ to a
representation of $\tilde{H}$ also denoted
by $\r_{\n,H}$.
We let $C = C_{\rm min}$ be the minimal
$H$-invariant 
cone
in $\fq$ generated by $X^0$. We denote by
$\tilde{C} = \tilde{C}_{\rm min}$ the minimal $G^c$-invariant
cone in $i\fg^c$. Then $\tilde{C}\cap
\fq = \pr_{\fq}(\tilde{C} = C$,
where $\pr_{\fq} : \fg \pf \fq$ denotes
the orthogonal projection (cf. \cite{HO95,'O90b}). As
$L_{\rm pos}\le 0$ it follows that
$\r_\l$ extends to a holomrophic representation
of the universal semigroup $\G (G^c,\tilde{C})$ corresponding
to $G^c$ and $\tilde{C}$, (cf. \cite{HiNe93}).
Let $G^c_1$ be the analytic subgroup of $G_{\bbC}$ corresponding
to the Lie algebra $\fg^c$. Let $H_1$ be the
analytic subgroup of $G^c_1$ corresponding to
$\fh$. Then - as we are assuming that $G\subset G_{\bbC}$ -
we have $H_1 = H_o$. Let
$\k : G^c \pf G^c_1$ be the canonical projection and
let 
$Z_H = \k^{-1}(Z_{G^c_1}\cap H_o)$. Then
$\r_{\n}$ is trivial on $Z_H$ as $\n - \r$ is trivial
on $\exp ([\fk^c,\fk^c])\supset H^c\cap K^c$. Thus
$\r_{\n}$ factor to
$G^c\bftdiagup Z_H$ and $\G (G^c,\tilde{C})\bftdiagup Z_H$.
Notice that $(G^c\bftdiagup Z_H)^\t_o$ is isomorphic to
$H_o$. Therefore we can view $H_o$ 
as subgroup of $G^c\bftdiagup Z_H$ and and $S_o(C) = H_o\exp C$
as a subsemigroup of $\G (G^c,\tilde{C})\bftdiagup Z_H$. In
particular
$\t_{\n}(s)$ is defined for $s\in S_o(C)$. This allows
us to write $\r_\n (h)$ or $\r_{\n,H}(h)$
for $h\in H_o$.
As $\fn_{\bbC} = \fp^+$ and $\fp^- = \halo{\fn}_{\bbC}
\t (\fn_{\bbC})$ it follows, using Lemma \ref{ConjugationCoincides}, that 
\[ a_{\halo{N}}(h)^{\n - \r} = \ip{u}{\r_{\n,H}(h)u}\, .\]

In particular we get that $(h,k)\mapsto a_{\halo{N}}(h^{-1}k)^{\n -\r}$
is positive semidefinite if $\n -\r \le L_{\rm pos}$.

\npar
Let
us now concider the case
$G = H_{\bbC}$ and $G^c = \tilde{H}\times \tilde{H}$.
Denote the constant $L_{\rm pos}$ for $\tilde{H}$ by
$S_{\rm pos}$ and denote for $\m \le S_{\rm pos}$ the
representation with lowest $\tilde{H}\cap \tilde{K}$-type
$\m$ by $(\t_\m,L_\m)$. Let $\halo{\t}_\m$ be the conjugate
representation. Recall that we view $\tilde{H}$ as a
subset of $G^c$ by the diagonal embedding
\[\tilde{H}\ni h\mapsto (h,h) \in \D (G^c) :=
\{(x,x)\in G^c\mid x\in \tilde{H}\}\, .\]
The center of $\fk^c$ is two dimensional (over $\bbR$)
and generated by $i(X^0,X^0)$ and $i(X^0,-X^0)$. We
choose $Z^0 = i(X^0,- X^0)$. Then
$\fp^+ = \fn\times \halo{\fn}$. Let
$u$ be again a lowest weight vector of norm one. 
Denote the corresponding vector in the conjugate Hilbert space
by $\halo{u}$. Then for $h\in \tilde{H}$:

\begin{eqnarray*}
\ip{u\otimes \halo{u}}{\t_\m\otimes \halo{\t}_\m(h,h)u\otimes \halo{u}}
&=& \ip{u}{\t_\l(h)u}\overline{ \ip{u}{\t_\l(h)u}}\\
&=& |\ip{u}{\t_\l(h)u}|^2\\
&=& a_{\halo{N}}(h)^{2\m}
\end{eqnarray*}

Thus we define in this case $L_{\rm  pos} := 2 S_{\rm pos}$.
As before we notice that $\t_\n\otimes\halo{\t}_\n (h,h)u\otimes
\halo{u}$ is well defined on $H$.
We how have:

\begin{Lemma}\label{ahaloN}
Assume that $\GH$ is non-compactly causal. For
$\n-\r\le L_{\rm pos}$ there exists an unitary irreducible
highest weight representation $(\r_\n,\bK_\n)$ of
$G^c$ and a lowest $K^c$-type vector $u$ of norm one such that
for every $h\in H$

\[a_{\halo{N}}(h)^{\n - \r} = \ip{u}{\t_\n(h)u}\, .\]
Hence the kernel
\[(H\times H)\ni (h,k)\mapsto a_{\halo{N}}(h)^{\n -\r}\in \bbR\]
is positive semidefinite. In
particular
$\ip{\cdot }{\cdot }_J$ is positive semidefinite on $\bC_c^\infty(\gO)$
for $\l - \r \le L_{\rm pos}$.
\end{Lemma}

\npar
The Basic Lemma and the L\"uscher-Mack
Theorem, together with the above, now imply the following Theorem

\begin{Satz}[Reflection Symmetry for Complementary Series]\label{PR}
Assume that $\GH$ is non-compactly causal
and such that there exists a $w\in K$ such
that $\Ad (w)|_{\fa} = -1$. Let
$\p_\n$ be a complementary series such that
$\n \le L_{\rm pos}$. Let $C$ be the minimal
$H$-invariant cone in $\fq$ such that $S(C) $ is
contained in the contraction
semigroup of $HP_{\rm max}$ in $G/P_{\rm max}$.
Let $\gO$ be the bounded realization of
$H\bftdiagup H\cap K$ in $\halo{\fn}$.
Let $J(f)(x) : = f(\t (x)w^{-1})$. Let
$\bK_0$ be the closure of $\cC_c^\infty(\gO)$ in
$\bH_\n$. Then the following holds:
\begin{enumerate}
\item[\hss\llap{\rm 1)}]  $(G,\t, \p_\n,C,J,\bK_0)$ satisfies
the positivity
conditions {\rm (PR1)--(PR2)}.
\item[\hss\llap{\rm 2)}] $\p_\n$ defines a contractive
representation $\tilde{\p}_\n $ of
$S(C)$ on $\bK$ such that $\tilde{\p}_\n(\g)^*
= \tilde{\p}_\n(\t (\g)^{-1})$.
\item[\hss\llap{\rm 3)}]  There exists a unitary representation
$\tilde{\p}^c_\n$ of $G^c$ such that
\begin{enumerate}
\item[\hss\llap{\rm i)}] $d\tilde{\p}_\n^c(X) = d\tilde{\p}_\n (X)\,
\quad \forall X\in \fh$.
\item[\hss\llap{\rm ii)}] $d\tilde{\p}_\l^c(iY) = i\, d\tilde{\p}_\l
(Y)\, \quad \forall Y\in C$.
\end{enumerate}
\end{enumerate}
\end{Satz}

\npar
We remark that this Theorem includes the results of
R. Schrader for $\lSL (2n,\bbC)\bftdiagup \SU (n,n)$, \cite{Sch86}.

\npar
We will now generalize this to all non-compactly
causal symmetric spaces and all $\n $ such
that $\n - \r \le L_{\rm pos}$. We will also
show that actually $\tilde{\p}_\n^c\simeq \r_\n$,
where $\r_\n$ is the irreducible unitary highest weight
representation of $G^c$ such that
\[a(h)^{\n - \r} = \ip{u}{\r_\n(h)u}\]
as before. {}From now on we assume that $\n - \r \le L_{\rm pos}$.
Let $\bK_0$ be the completion of $\cC_c^\infty (\gO)$ in
the norm $\ip{\cdot}{A(\n)J(\cdot )}$. Let
$\bN$ be the space of vectors of zero lenght and
let $\bK $ be the completion of $\bK_0/\bN$ in the
induced norm. First of all we have to show that
$\p_\n(\g)$ passes to a continuous operator
$\tilde{\p}_\n(\g)$ on $\bK$ such
that $\tilde{\p}_\n(\g)^* = \tilde{\p}_\n(\t (\g)^{-1})$.
For that we recall that 
\begin{equation}
H\bftdiagup H\cap K
= H_o\bftdiagup H_o\cap K = \gO
\end{equation}
so we my replace the integration over $H$ in $\ip{f}{A(\n)Jf}_\n$
with integration over $H_o$. For $f\in \cC^\infty_c(\gO)$ define
\begin{equation}
\r_\n(f)u := \int_{H_o} f(h\cdot 0)\r_\n(h)u\, dh\, .
\end{equation}

\begin{Lemma}\label{L:posref} Assume that $\n - \r \le L_{\rm pos}$.
Let $\r_\n,\, \bK_\n$ and $u$ be as specified in Lemma \ref{ahaloN}.
and let $f, g\in \cC^\infty_c(\gO)$ and $s\in S(C)$.
Then the following holds:
\begin{enumerate}
\item[\hss\llap{\rm 1)}] $\ip{f}{[A(\n)J](g)}_\n
= \ip{\r_\n(f)u}{\r_\n(g)u}$.
\item[\hss\llap{\rm 2)}] $\r_\n(\p_\n(s)f)u = \r_\n(s)\r_\n(f)u$.
\item[\hss\llap{\rm 3)}] $\p_\n(s)$ passes to a contractive operator
$\tilde{\p}_\n(s)$ on $\bK$ such that $\tilde{\p}_\n(s)^*
= \tilde{\p}_\n(\t (s)^{-1})$.
\end{enumerate}
\end{Lemma} 

{\it Proof\/}: (1) Let $f$ and $g$ be as above. Then
\begin{eqnarray*}
\ip{f}
{[A(\n)J](g)} &=& \int_{H_o\bftdiagup H_o\cap K}\int_{H_o\bftdiagup H_o\cap K}
\overline{f(h)}g(k)a_{\halo{N}}(h^{-1}k)^{\n-\r}\, dh\, dk\\
&=& \int_{H_o\bftdiagup H_o\cap K}\int_{H_o\bftdiagup H_o\cap K}
\overline{f(h)}g(k)\ip{u}{\r_\n(h^{-1}k)u}\, dh\, dk\\
&=& \int_{H_o\bftdiagup H_o\cap K}\int_{H_o\bftdiagup H_o\cap K}
\overline{f(h)}g(k)\ip{\r_\n(h)u}{\r_\n(k)u}\, dh\, dk\\
&=& \ip{\r_\n(f)u}{\r_\n(g)}
\end{eqnarray*}
This proves (1).
\npar
(2) This follows from Lemma \ref{L:Int},7) and the following
calculation:
\begin{eqnarray*}
\r_\n(\p_\n(s)f)u &=& \int f(s^{-1}h)\r_\n(h)u\, dh\\
&=& \int f(h(s^{-1}h))a_H(s^{-1}h)^{-(\n + \r)}
\r_\n(h)u\, dh\\
&=& \int f(h(s^{-1}h))a_H(sh(s^{-1}))^{\n -\r}\r_\n(h)
a_H(s^{-1}h)^{-2\r} u \, dh\\
&=& \int f(h)\r_\n(sh)u\, dh\\
&=& \r_\n(s)\r_\n(f)u\, ,
\end{eqnarray*}
where we have used that
\[\r_\n(sh)u = a_H(sh)^{\n -\r}\r_\n(h(sh))u\, .\]
\npar
(3) By (1) and (2) we get: 
\begin{eqnarray*}
\|\p_\n(s)f\|_J^2 & = &\|\r_\n (s)\r_\n(f)u\|^2 \le
\|\r_n(f)u\|^2\\
& = & \ip{f}{[A(\n )J]f}_\n (= \|f\|_J^2)
\end{eqnarray*}
Thus $\p_\n (s)$ passes to a contractive operator on
$\bK$. That $\tilde{\p}_\n(s)^*
= \tilde{\p}_\n(\t (s)^{-1})$ follows from
Lemma \ref{L:AnuJ}.
\hfill$\Box$\medskip

\begin{Satz}[Identification Theorem]\label{S:Posref} Assume
that $\GH$ is non-compactly
causal and that $\n - \r \le L_{\rm pos}$.
Let $\r_\n$, $\bK_\n$ and $u\in \bK_\n$ be as in Lemma
\ref{ahaloN}.
Then the following hold:

\begin{enumerate}
\item[\hss\llap{\rm 1)}] There exists a continuous
contractive representation
$\tilde{\p}_\n$ of $S_o(C)$ on $\bK$ such that
\[ \tilde{\p}_\n(s)^*
= \tilde{\p}_\n(\t (s)^{-1})\, ,\quad \forall s\in S_o(C)\, .\]
\item[\hss\llap{\rm 2)}] There exists a unitary representation
$\tilde{\p}^c_\n$ of $G^c$ such that
\begin{enumerate}
\item[\hss\llap{\rm i)}] $d\tilde{\p}_\n^c(X) = d\tilde{\p}_\n (X)\,
\quad \forall X\in \fh$.
\item[\hss\llap{\rm ii)}] $d\tilde{\p}_\n^c(iY) = i\, d\tilde{\p}_\n
(Y)\, \quad \forall Y\in C$.
\end{enumerate}
\item[\hss\llap{\rm 3)}] The
map
\[\cC^\infty_c(\gO)  \ni f \mapsto \r_\n(f)u \in \bK_\n\]
extends to an isometry $\bK \simeq \bK_\n$
intertwining $\tilde{\p}^c_\n$ and $\r_\n$.
In particular $\tilde{\p}^c_\n$ is irreducible
and isomorphic to $\r_\n$.
\end{enumerate}
\end{Satz}

{\it Proof\/}:  (1)  follows from Lemma \ref{L:posref}
as obviously $\tilde{\p}_\n (sr) = \tilde{\p}_\n (s)\tilde{\p}_\n(r)$.
\npar
(2) This follows now from the Theorem of L\"uscher-Mack.
\npar
(3) By Lemma \ref{L:posref} we know
that $f\mapsto \r_\n(f)u$ defines an isometric
$S_o(C)$-intertwining operator. Let $f\in \cC_c^\infty(\gO)$.
Differentiation and the fact that
$\t_\n$ is holomorphic gives

\begin{enumerate}
\item[\hss\llap{\rm i)}] $\r_\n(d\tilde{\p}_\n^c(X)f)u =
d\r_\n (X)\r_\n(f)u\, ,
\quad \forall X\in \fh$.
\item[\hss\llap{\rm ii)}] $\r_\n(i\, d\tilde{\p}_\n^c(Y)f)u =
i\, d\r_\n(Y)\r_\n(f)u\, ,\quad \forall Y\in C$.
\end{enumerate}
But those are exactly the relations that define $\tilde{\p}_\n^c$.
The fact that $\fh\oplus iC$ generates $\fg^c$ implies that
$f\mapsto \r_\n(f)u$ induces an $\fg^c$-intertwining
operator
intertwining $\tilde{\p}^c_\n$ and $\r_\n$. As both are
also representations of $G^c$, it follows that this
is an isometric $G^c$-map. In particular if this is not
the zero-map it has to be an isomorphism as
$\r_\n$ is irreducible.
Choose a sequence $\{f_j\}$ in $\cC^\infty_c(\gO)$
approximating the Delta function. The usual calculation shows that
\[\r_\n(f_j)u\to u\, .\]
Hence there is a $j$ such that $\r_\n(f_j)u\not= 0$. This
proves the theorem. \hfill$\Box$\medskip

\begin{Rem}{\rm This above Theorem  realise the highest
weight representation $\r_\n$ on a function space on
$H/H\cap K$. The construction is in some sense inverse to
the construction in \cite{OO96}. The highest weight representation
$\r_\n$ can be realized in a Hilbert space $\cO$ of holomorphic
functions on $\gO_{\bbC}$. The restriction of
a holomorphic function to $\gO$ is injective by Lemma \ref{T:compdom}.
Multiplying by a suitable character induces then
a injective $H$-intertwining operator into $\bL^2(H\bftdiagup H\cap K)$,
at least for $\n$ big enough. We refer to \cite{OO96}
for further details.
}
\end{Rem}
 \npar
We will now explain another view of the above results using
local representations instead on the Theorem of
L\"uscher-Mack. This will use the realization
as functions on $\gO\subset \halo{\fn}$ and
in particular explain the
kernel 
\[ (X,Y)\mapsto K_\n (X,Y) := a_{\halo{N}}(\t (\exp X )^{-1}\exp Y )^{\n -\r}
\, \quad X,Y\in \gO\, .\]

For this we recall some results from \cite{MDRF96}, in particular
Theorem 5.1 and Theorem 7.1:
We assume that $\n-\r\le L_{\rm pos}$. Let
$\r_\n$, $\bK_\n$ and $u$ be as before. Then
\[ K_\n(X,Y) = \ip{u}{\r_\n(\exp (-\t (X))\exp (Y))u}
= \ip{\r_\n(\exp (X))u}{\r_\n(\exp (Y))u}\]
because of Lemma \ref{ConjugationCoincides}. 

\begin{Lemma}\label{ReproducingKernel} Let the notation be
as above. Then the following holds.
\begin{enumerate}
\item[\hss\llap{\rm 1)}] The map
\[\gO\ni X \mapsto q_X u :=  \r_\n (\exp X)u\in \bK_\n \]
extends to a holomorphic map
on $\gO_{\bbC}$ given by
\[q_X u = \sum_{n=0}^\infty \frac{d\r_\n(X)^nu}{n!}\, .\]
\item[\hss\llap{\rm 2)}]  The function
$\ip{ q_X u }{q_Yu}$ is 
an extension of $K_\n(X,Y)$ to $\gO_{\bbC}\times \gO_{\bbC}$,
holomorphic in the second variable and antiholomorphic in the
first variable. We will denote this extension also by
$K_\n(X,Y)$.
\item[\hss\llap{\rm 3)}] The function $K_\n(X,Y)$ is positive
definite.
\end{enumerate}
\end{Lemma}

{\it Proof\/}: See \cite{MDRF96}.\hfill$\Box$\medskip

\npar
Let
$U \subset \gO$ be open.
We identify $\cC^\infty_c(U)$ with the space
of elements in $f\in \cC^\infty(\gO)$ such that
$f|_{\halo{N}}\circ\z^{-1}\in \cC^\infty_c(U)$. For
$R>0$,  let
\[B_R := \Ad (H\cap K) \{\sum_{j=1}^r t_jX_{-j}\mid
-R < t_j < R\}\, .\] Then $B_R$ is open in $\halo{\fn}$.  Let $\b :\bK
_0\pf \left( \bK _0\bftdiagup \bN \right) \tilde{}=\bK $ be the
canonical map. Then $\b $ is a contraction ($\left\| \b (f)\right\|
^2_J=\ip{f}{Jf}\le \left\| f\right\| ^2$). For $U\subset \gO $ open,
let $\bK (U):=\b (\cC _c^\infty (U))$.
\begin{Satz}\label{Dense}Let $U\subset \gO$ be open. Then
$\bK (U)$ is dense in $\bK $.
\end{Satz}

{\it Proof\/}: Let $x\in U$. Then we can choose $h\in H$
such that $hx = 0$. As
$\cC^\infty_c(U) = h\cdot C^\infty_c(h \cdot U)$ and
$H$ acts unitarily, it follows that we can assume
that $0\in U$. Let $R>0$ be such that
$B_R\subset U$. Then $\cC^\infty_c(B_R)\subset \cC^\infty_c(U)$.
Hence we can assume that $U = B_R$.
Let
$g\in \cC^\infty_c(U)^\perp$ and let $f\in \cC^\infty_c(\gO)$.
We want to show that $\ip{g}{f}_J = 0$.
Choose $0<L<1$ such that $\Supp (f)\subset B_L$. For
$t\in \bbR $ and $a_t = \exp (2tX^0)$ we have
$a_t\cdot B_L = B_{e^{-2t}L}$. Thus
$\Supp (\p (\n)(a_t)f) \subset B_{e^{-2t}L}$. Choose
$0<s_0$ such that $e^{-2t}L < R$ for every $t> s_0$.
Then $\p (\n)(a_t)(f)\in \cC^\infty_c(U)$ for every
$t>s_0$. It follows that for $t>s_0$:
\begin{eqnarray*}
0 &=& \ip{g}{\p (\n )(a_t)f}_J\\
&=& \int \int \overline{g(x)}\left[ \p (\n )(a_t)f\right] (y)K_\l
(x,y)\,dx\,dy\\
&=& e^{(\l +1)t}\int \int \overline{g(x)}f(e^{2t}y)K_\l (x,y)\,dx\,dy\\
&=& e^{(\l -1)t}\int \int \overline{g(x)}f(y)K_\l (x,e^{-2t}y)\,dx\,dy\, .
\end{eqnarray*}
By Lemma \ref{ReproducingKernel} we know that $z\mapsto K_\l (x,zy)$
is holomorphic on $D=\left\{ z\mid \left| z\right|
<1\right\} $. As $g$ and $f$ both have compact support it follows that
\[F(z):=\int \int \overline{g(x)}f(y)K_\l (x,zy)\,dx\,dy\]
is holomorphic on $D$. But $F(z)=0$ for $0<z<e^{-2s_0}$. Thus $F(z)=0$
for every $z$. In particular
\[\ip{g}{\p (\n )(a_t)f}_J=0\]
for every $t>0$. By continuity $\ip{g}{f}_J=0$. Thus $g=0$.\hfill$\Box$\medskip

\npar
Let us recall some basic facts from \cite{Jor86}. Let 
$\r$ be a local homomorphism of a neighborhood
$U$ of $e$ in $G$ into 
the space of linear operators on the
Hilbert space $\bH$ such that $\r (g)$ is densely defined
for $g\in U$. Furthermore $\r |_{(U\cap H)}$ extends
to a strongly continuous representation of $H$ in
$\bH$.
$\r $ is called a {\em local representation} if there
exists a dense subspace $\bD\subset \bH$ such that 
the following holds:
\begin{list}{}{\setlength{\leftmargin}{\customleftmargin} 
\setlength{\itemsep}{0.5\customskipamount}
\setlength{\parsep}{0.5\customskipamount}
\setlength{\topsep}{\customskipamount}}
\item[\hss\llap{\rm LR1)}] $\forall g\in U$, $\bD \subset \bD (\r (g))$, where
$\bD (\r (g))$ is the domain of definition for $\r (g)$.
\item[\hss\llap{\rm LR2)}] If $g_1,g_2,g_1g_2\in U$ and $u\in \bD$ then
$\r (g_2)u \in \bD (\r (g_1))$ and
\[\r (g_1)[\r (g_2)u] = \r (g_1g_2)u\, .\]
\item[\hss\llap{\rm LR3)}] Let $Y\in \fh$ such that $\exp tY\in U$ for
$0\le t\le 1$. Then for every $u\in \bD$ 
\[\lim_{t\to 0} \r (\exp tY) u = u\, .\]
\item[\hss\llap{\rm LR4)}] $\r (Y)\bD\subset\bD$ for every $Y\in \fh$.
\item[\hss\llap{\rm LR5)}] $\forall u\in \bD\, \exists V_u$ an
open $1$-neighborhod in $H$ such that $UV_u\subset U^2$ and
$\r (h)u \in \bD$ for every $h\in V_u$.
\item[\hss\llap{\rm LR6)}] For every $Y\in \fq$ and every $u\in \bD$ the
function
\[h\mapsto \r (\exp (\Ad (h)Y)u\]
is locally integrable on $\{h\in H\mid \exp (\Ad (h)Y)\in U\}$.
\end{list}
\cite{Jor86} now states that every local representation
extends to a unitary representation of $G^c$.  We now want to use
Theorem \ref{Dense} to construct a local representation of $G$. For
that let $0<R<1$ and let $\bD =\bK (B_R(0))$. Let $V$ be a symmetric
open neighborhood of $1\in G$ such that $V\cdot B_R(0)\subset \gO
$. Let $U_1$ be a convex symmetric neighborhood of $0$ in $\fg $ such
that with $U:=\exp U_1$ we have $U^2\subset V$. If $g\in U$ then
obviously (LR1)--(LR3) are satisfied. (LR4) is satisfied as
differentiation does not increase support. (LR6) is also clear as
$u=\b (f)$ with $f\in \cC _c^\infty (U)$ and hence $\left\| \r _c(\exp
\Ad (h)Y)u\right\| $ is continuous as a function of $h$.

\npar (LR5) Let $u=\b (f)\in \bK (B_R(0))$. Let $L=\Supp (f)\subset
B_R(0)$. Let $V_u$ be such that $V_u^{-1}=V_u$, $V_uL\subset B_R(0)$,
and $V_u\subset U$. Then $UV_u\subset U^2$ and $\tilde{\p }(\n
)(h)u=\b (\p (\n )(h)f)$ is defined and in $D$. This now implies that
$\tilde{\p }$ restricted to $U$ is a local representation. Hence the
existence of $\tilde{\p }^c$ follows from \cite{Jor86}. We notice that
this construction of $\tilde{\p }^c$ does not use the full semigroup
$S$ but only $H$ and $\exp \bbR _+X^o$.\hfill$\Box$\medskip

\section{\protect\label{Diagonal}The diagonal case $\p\oplus (\p\circ \t)$}
\setcounter{equation}{0}

A special case of the setup in Definition \ref{ReflectionSymmetric}
above arises as follows: Let the group $G$, and $\t \in \Aut _2(G)$ be
as described there. Let $\bH _\pm $ be two given complex Hilbert
spaces, and $\p _\pm \in \mathop{{\rm Rep}} (G,\bH _\pm )$ a pair of
unitary representations. Suppose $T:\bH _-\pf \bH _+$ is a unitary
operator such that $T\p _-=\left( \p _+\circ \t \right) T$, or
equivalently,
\begin{equation}\label{TPiMinus}
T\p _-(g)f_-=\p _+\left( \t (g)\right) Tf_-
\end{equation}
for all $g\in G$, and all $f_-\in \bH _-$. Form the direct sum $\bH
:=\bH _+\oplus \bH _-$ with inner product
\begin{equation}\label{IPPlusMinus}
\ip{f_+\oplus f_-}{f^{\prime }_+\oplus f^{\prime
}_-}:=\ip{f_+}{f^{\prime }_+}_++\ip{f_-}{f^{\prime }_-}_-
\end{equation}
where the $\pm $ subscripts are put in to refer to the respective
Hilbert spaces $\bH _\pm $, and we may form $\p :=\p _+\oplus \p _-$
as a unitary representation on $\bH =\bH _+\oplus \bH _-$ by
\[\p (g)\left( f_+\oplus f_-\right) =\p _+(g)f_+\oplus \p _-(g)f_-\,
,\quad g\in G,\; f_\pm \in \bH _\pm \, .\]
Setting
\begin{equation}\label{JMatrix}
J:=\left( \matrix{ 0 & T\cr T^* & 0\cr }\right) \, ,
\end{equation}
i.e., $J\left( f_+\oplus f_-\right) =\left( Tf_-\right) \oplus \left(
T^*f_+\right) $, it is then clear that properties (1)--(2) from
Definition \ref{ReflectionSymmetric} will be satisfied for the pair
$(J,\p )$. Formula (\ref{TPiMinus}) may be recovered by writing out
the relation
\begin{equation}\label{JPiRelation}
J\p =\left( \p \circ \t \right) J
\end{equation}
in matrix form, specifically
\[\left( \matrix{ 0 & T\cr T^* & 0\cr }\right) \left( \matrix{
\p _+(g) & 0\cr 0 & \p _-(g)\cr }\right) =\left( \matrix{ \p _+(\t
(g)) & 0\cr 0 & \p _-(\t (g))\cr }\right) \left( \matrix{ 0 & T\cr
T^* & 0\cr }\right) \, .\]
If, conversely, (\ref{JPiRelation}) is assumed for some unitary
period-$2$ operator $J$ on $\bH =\bH _+\oplus \bH _-$, and, if the two
representations $\p _+$ and $\p _-$ are {\it disjoint,} in the sense
that no irreducible in one occurs in the other (or, equivalently,
there is no nonzero intertwiner between them), then, in fact,
(\ref{TPiMinus}) will follow from (\ref{JPiRelation}). The diagonal
terms in (\ref{JMatrix}) will be zero if (\ref{JPiRelation})
holds. This last implication is an application of Schur's lemma.

\begin{Lemma}\label{Operator}Let $0\ne \bK _0$ be a closed linear
subspace of $\bH =\bH _+\oplus \bH _-$ satisfying the positivity
condition {\rm (PR3)} in Definition \ref{Hyperbolic}, i.e.,
\begin{equation}\label{vJvPositivity}
\ip{v}{Jv}\ge 0\, ,\quad \forall v\in \bK _0
\end{equation}
where
\begin{equation}\label{JMatrixBis}
J=\left( \matrix{ 0 & T\cr T^* & 0\cr }\right) 
\end{equation}
is given from a fixed unitary isomorphism $T:\bH _-\pf \bH _+$ as in
{\rm (\ref{TPiMinus}).} For $v=f_+\oplus f_- \in \bH =\bH _+\oplus \bH
_-$, set $P_+v:=f_+$. The closure of the subspace $P_+\bK _0$ in $\bH
_+$ will be denoted $\overline{P_+\bK _0}$. Then the subspace
\[\bG =\left\{ \left. \left( \matrix{ f_+\cr f_-\cr }\right) \in
\bK _0\,\right|\, f_-\in T^*\left( \overline{P_+\bK _0}\right) \right\}\]
is the graph of a closed linear operator $M$ with domain
\begin{equation}\label{MDomain}
\bD =\left\{ f_+\in \bH _+\,\left| \,\exists f_-\in T^*\left(
\overline{P_+\bK _0}\right) \mathop{{\rm s.t.}}\left( \matrix{
f_+\cr f_-\cr }\right) \in \bK _0\right. \right\}\, ;
\end{equation}
and, moreover, the product operator $L:=TM$ is dissipative on this
domain, i.e.,
\begin{equation}\label{LDissipative}
\ip{Lf_+}{f_+}_++\ip{f_+}{Lf_+}_+\ge 0
\end{equation}
holds for all $f_+\in \bD $.
\end{Lemma}

{\it Proof\/}: The details will only be sketched here, but the reader
is referred to \cite{Sto51} and \cite{Jor80} for definitions and
background literature. An important argument in the proof is the
verification that, if a column vector of the form $\left( \matrix{
0\cr f_-\cr }\right) $ is in $\bG $, then $f_-$ must necessarily be
zero in $\bH _-$. But using positivity, we have
\begin{equation}\label{uJvBound}
\left| \ip{u}{Jv}\right| ^2\le \ip{u}{Ju}\ip{v}{Jv}\, ,\quad \forall
u,v\in \bK _0\, .
\end{equation}
Using this on the vectors $u=\left( \matrix{ 0\cr f_-\cr }\right) $
and $v=\left( \matrix{ k_+\cr k_-\cr }\right) \in \bK _0$, we get
\[\ip{\left( \matrix{ 0\cr f_-\cr }\right) }{\left( \matrix{ Tk_-\cr
T^*k_+\cr }\right) }=\ip{f_-}{T^*k_+}=0\, ,\quad \forall k_+=P_+v\,
.\] But, since $f_-$ is also in $T^*\left( \overline{P_+\bK _0}\right)
$, we conclude that $f_-=0$, proving that $\bG $ is the graph of an
operator $M$ as specified. The dissipativity of the operator $L=TM$ is
just a restatement of {\rm (PR3)}.\hfill$\Box$\medskip

\npar The above result involves only the operator-theoretic
information implied by the data in Definition \ref{Hyperbolic}, and,
in the next lemma, we introduce the representations:

\begin{Lemma}\label{Normal}Let the representations $\p _\pm $ and
the intertwiner $T$ be given as specified before the statement of
Lemma \ref{Operator}. Let $H=G^\t $; and suppose we have a 
cone $C\subset \fq$ as specified in {\rm (PR2${}^\prime$)}.
Assume {\rm (PR1), (PR2${}^\prime$) and (PR3${}^\prime$)} and
assume further that
\begin{list}{}{\setlength{\leftmargin}{\customleftmargin} 
\setlength{\itemsep}{0.5\customskipamount}
\setlength{\parsep}{0.5\customskipamount}
\setlength{\topsep}{\customskipamount}}
\item[\hss\llap{{\rm PR4)}}]  $\bD $ is dense in $\bH _+$, and
\item[\hss\llap{{\rm PR5)}}]  the commutant of $\left\{ \p _+(h)\mid
h\in H\right\} $ is abelian.
\end{list}
Then it follows that the operator $L=TM$ is normal.
\end{Lemma}

{\it Proof\/}: Since $T$ is a unitary isomorphism $\bH _-\pf \bH _+$
we may make an identification and reduce the proof to the case where
$\bH _+=\bH _-$ and $T$ is the identity operator. We then have
\[\p _-=T^{-1}\left( \p _+\circ \t \right) T=\p _+\circ \t \, ;\]
and if $h\in H$, then
\[\p _-(h)=\p _+\left( \t (h)\right) =\p _+(h)\, ;\]
while, if $\t (g)=g^{-1}$, then
\[\p _-(g)=\p _+\left( \t (g)\right) =\p _+\left( g^{-1}\right) \, .\]
Using only the $H$ part from (PR2${}^\prime$), we conclude
that $\bK _0$ is invariant under 
$\p _+\oplus \p _+ (H)$. If the projection $P_{\bK _0}$ of $\bH _+\oplus \bH
_+$ onto $\bK _0$ is written as an operator matrix $\left( \matrix{
P_{11} & P_{12}\cr P_{21} & P_{22}\cr }\right) $ with entries
representing operators in $\bH _+$, and satisfying
\begin{eqnarray*}
P_{11}^* &=& P_{11}\, ,\\
P_{22}^* &=& P_{22}\, ,\\
P_{12}^* &=& P_{21}\, ,\\
P_{ij} &=& P_{i1}P_{1j}+P_{i2}P_{2j}\, ,
\end{eqnarray*}
then it follows that
\begin{equation}\label{PijPiPlus}
P_{ij}\p _+(h)=\p _+(h)P_{ij}\quad \forall i,j=1,2,\;\forall h\in H\, ,
\end{equation}
which puts each of the four operators $P_{ij}$ in the commutant $\p
_+(H)^{\prime }$ from (PR5). Using (PR4), we then conclude that $L$ is
a dissipative operator with $\bD $ as dense domain, and that $\bK _0$
is the graph of this operator. Using (PR5), and a theorem of Stone
\cite{Sto51}, we finally conclude that $L$ is a normal operator, i.e.,
it can be represented as a multiplication operator with dense domain
$\bD $ in $\bH _+$.\hfill$\Box$\medskip

\npar We shall consider two cases below (the Heisenberg group, and the
$(ax+b)$-group) when conditions (PR4)--(PR5) can be verified from the
context of the representations. Suppose $G$ has two abelian subgroups
$H$, $N$, and the second $N$ also a normal subgroup, such that $G=HN$
is a product representation in the sense of Mackey \cite{Mac}. Define
$\t \in \Aut _2(G)$ by setting
\begin{equation}\label{TauDef}
\t (h)=h\, ,\quad \forall h\in H\, ,\und \t (n)=n^{-1}\, ,\quad
\forall n\in N\, .
\end{equation}

\npar The Heisenberg group is a copy of $\bbR ^3$ represented as
matrices $\vphantom{\left( \matrix{ 1 & a & c\cr 0 & 1 & b\cr 0 & 0 &
1\cr }\right) _X}\left( \matrix{ 1 & a & c\cr 0 & 1 & b\cr 0 & 0 &
1\cr }\right) $, or equivalently vectors $(a,b,c)\in \bbR ^3$. Setting
$H=\left\{ (a,0,0)\mid a\in \bbR \right\} $ and
\begin{equation}\label{NHeisenberg}
N=\left\{ (0,b,c)\mid b,c\in \bbR \right\} \, ,
\end{equation}
we arrive at one example.

\npar The $(ax+b)$-group is a copy of $\bbR ^2$ represented as
matrices $\left( \matrix{ a & b\cr 0 & 1\cr }\right) $, $a=e^s$, $b\in
\bbR $, $s\in \bbR $. Here we may take $H=\left\{ \left. \left(
\matrix{ a & 0\cr 0 & 1\cr }\right) \,\right| \,a\in \bbR _+\right\} $
and
\begin{equation}\label{Naxb}
N=\left\{ \left. \left( \matrix{ 1 & b\cr 0 & 1\cr }\right) \,\right|
\,b\in \bbR \right\} \, ,
\end{equation}
and we have a second example of the Mackey factorization. Generally,
if $G=HN$ is specified as described, we use the representations of $G$
which are induced from one-dimensional representations of $N$. If $G$
is the Heisenberg group, or the $(ax+b)$-group, we get all the
infinite-dimensional irreducible representations of $G$ by this
induction (up to unitary equivalence, of course). For the Heisenberg
group, the representations are indexed by $\hslash \in \bbR
\bftdiagdown \{ 0\} $, $\hslash $ denoting Planck's constant. The
representation $\p _\hslash $ may be given in $\bH = \bL ^2(\bbR)$ by
\begin{equation}\label{PiHBar}
\p _\hslash (a,b,c)f(x)=e^{i\hslash (c+bx)}f(x+a)\, ,\quad \forall
f\in \bL ^2(\bbR ),\; (a,b,c)\in G\, .
\end{equation}
The Stone-von Neumann uniqueness theorem asserts that every unitary
representation $\p $ of $G$ satisfying
\[\p (0,0,c)=e^{i\hslash c}I_{\bH (\p )}\quad (\hslash \ne 0)\]
is unitarily equivalent to a direct sum of copies of the
representation $\p _\hslash $ in (\ref{PiHBar}).

\npar The $(ax+b)$-group (in the form $\left\{ \left. \left( \matrix{
e^s & b\cr 0 & 1\cr }\right) \,\right| \,s,b\in \bbR \right\} $) has
only two inequivalent unitary irreducible representations, and they
may also be given in the same Hilbert space $\bL ^2(\bbR )$ by
\begin{equation}\label{PiPlusMinusaxb}
\p _\pm \left( \matrix{ e^s & b\cr 0 & 1\cr }\right) f(x)=e^{\pm
ie^xb}f(x+s)\, ,\quad \forall f\in \bL ^2(\bbR )\, .
\end{equation}
There are many references for these standard facts from representation
theory; see, e.g., \cite{Jor88}.

\begin{Lemma}\label{Abel}Let the group $G$ have the form $G=HN$
for locally compact abelian subgroups $H,N$, with $N$ normal, and
$H\cap N=\{ e\} $. Let $\chi $ be a one-dimensional unitary
representation of $N$, and let $\p =\mathop{{\rm ind}}_N^G(\chi )$ be
the corresponding induced representation. Then the commutant of
$\left\{ \p (H)\mid h\in H\right\} $ is an abelian von Neumann
algebra: in other words, condition {\rm (PR5)} in Lemma \ref{Normal}
is satisfied.
\end{Lemma}

{\it Proof\/}: See, e.g., \cite{Jor88}.\hfill$\Box$\medskip

\npar In the rest of the present section, we will treat the case of
the {\em Heisenberg group,} and the $(ax+b)${\em -group} will be the
subject of the next section.

\npar For both groups we get pairs of unitary representations $\p _\pm
$ arising from some $\t \in \Aut _2(G)$ and described as in
(\ref{JPiRelation}) above. But when the two representations $\p _+$
and $\p _-=\p _+\circ \t $ are irreducible and disjoint, we will show
that there are no spaces $\bK _0$ satisfying (PR1), (PR2$^{\prime }$),
and (PR3) such that $\bK =\left( \bK _0\bftdiagup \bN \right)
\tilde{}$ is nontrivial. Here (PR2) is replaced by
\begin{list}{}{\setlength{\leftmargin}{\customleftmargin} 
\setlength{\itemsep}{0.5\customskipamount}
\setlength{\parsep}{0.5\customskipamount}
\setlength{\topsep}{\customskipamount}}
\item[\hss\llap{\rm PR2$^{\prime }$)}]  $C$ is a nontrivial cone in $\fq $.
\end{list}
Since for both groups, and common to all the representations, we noted
that the Hilbert space $\bH _+$ may be taken as $\bL ^2(\bbR )$, we
can have $J$ from (\ref{JMatrixBis}) represented in the form $J=\left(
\matrix{ 0 & I\cr I & 0\cr }\right) $. Then the $J$-inner product on
$\bH _+\oplus \bH _-=\bL ^2(\bbR )\oplus \bL ^2(\bbR ) \simeq \bL
^2(\bbR ,\bbC ^2)$ may be brought into the form
\begin{equation}\label{JInnerProduct} \ip{\left( \matrix{ f_+\cr
f_-\cr }\right) }{\left( \matrix{ f_+\cr f_-\cr }\right) }_J=2\Re
\ip{f_+}{f_-}=2\int _{-\infty }^\infty \Re \left(
\overline{f_+(x)}f_-(x)\right) \,dx\, .
\end{equation}

\npar
For the two examples, we introduce
\[N_+=\left\{ (0,b,c)\mid b,c\in \bbR _+\right\} \]
where $N$ is defined in (\ref{NHeisenberg}), but $N_+$ is not
$H$-invariant. Alternatively, set
\[N_+=\left\{ \left. \left( \matrix{ 1 & b\cr 0 & 1\cr }\right)
\,\right|\, b\in \bbR _+\right\} \] for the alternative case where $N$
is defined from (\ref{Naxb}), and note that this $N_+$ is
$H$-invariant. In fact there are the following $4$ invariant cones in
$\fq $:
\begin{eqnarray*}
C_1^+ &=& \left\{ (0,0,t)\mid t\ge 0\right\} \\
C_1^- &=& \left\{ (0,0,t)\mid t\le 0\right\} \\
C_2^+ &=& \left\{ (0,x,y)\mid x\in \bbR ,\;y\ge 0\right\} \\
C_2^- &=& \left\{ (0,x,y)\mid x\in \bbR ,\;y\le 0\right\} \\
\end{eqnarray*}
Let $\p $ denote one of the representations of $G=HN$ from the
discussion above (see formulas (\ref{PiHBar}) and
(\ref{PiPlusMinusaxb})) and let $\cD $ be a closed subspace of $\bH
=\bL ^2(\bbR )$ which is assumed invariant under $\p (HN_+)$. Then it
follows that the two spaces
\begin{eqnarray}
\cD _\infty :=\bigvee \left\{ \p (n)\cD \mid n\in N\right\} \label{DInfinity}\\
\cD _{-\infty }:=\bigwedge \left\{ \p (n)\cD \mid n\in N\right\}
\label{DMinusInfinity}
\end{eqnarray}
are invariant under $\p (G)$, where the symbols $\bigvee $ and
$\bigwedge $ are used for the usual lattice operations on closed
subspaces in $\bH $. We leave the easy verification to the reader, but
the issue is resumed in the next section. If $P_\infty $, resp.,
$P_{-\infty }$, denotes the projection of $\bH $ onto $\cD _\infty $,
resp., $\cD _{-\infty }$, then we assert that both projections $P_{\pm
\infty }$ are in the commutant of $\p (G)$. So, if $\p $ is
irreducible, then each $P_\infty $, or $P_{-\infty }$, must be $0$ or
$I$. Since $\cD _{-\infty }\subset \cD \subset \cD_\infty $ from the
assumption, it follows that $P_\infty =I$ if $\cD \ne \{ 0\} $.

\begin{Lemma}\label{Heisenberg}Let $G$ be the Heisenberg group,
and let the notation be as described above. Let $\p _+$ be one of the
representations $\p _\hslash $ and let $\p _-$ be the corresponding
$\p _{-\hslash }$ representation. Let $0\ne \bK _0\subset \bL ^2(\bbR
)\oplus \bL ^2(\bbR )$ be a closed subspace which is invariant under
$\left( \p _+\oplus \p _-\right) \left( HN_+\right) $. Then it follows
that there are only the following possiblities for $\overline{P_+\bK
_0}$: $\{ 0\}$, $\bL ^2(\bbR )$, or $A\cH _+$ where $\cH _+$ denotes
the Hardy space in $\bL ^2(\bbR )$ consisting of functions $f$ with
Fourier transform $\hat{f}$ supported in the half-line $\left[
0,\infty \right) $, and where $A\in \bL ^\infty (\bbR )$ is such that
$\left| A(x)\right| =1\mathop{{\rm a.e.}}x\in \bbR $. For the space
$\overline{P_-\bK _0}$, there are the possibilities: $\{ 0\}$, $\bL
^2(\bbR )$, and $A\cH _-$, where $A$ is a (possibly different) unitary
$\bL ^\infty $-function, and $\cH _-$ denotes the negative Hardy
space.
\end{Lemma}

{\it Proof\/}: Immediate from the discussion, and the Beurling-Lax
theorem classifying the closed subspaces in $\bL ^2(\bbR )$ which are
invariant under the multiplication operators, $f(x)\mapsto
e^{iax}f(x)$, $a\in \bbR _+$. We refer to \cite{LaPh}, or
\cite{Hel64}, for a review of the Beurling-Lax
theorem.\hfill$\Box$\medskip

\begin{Cor}\label{Positive}Let $\p _\pm $ be the representations
of the Heisenberg group, and suppose that the subspace $\bK _0$ from
Lemma \ref{Heisenberg} is chosen such that {\rm (PR1)--(PR3)} in
Definition \ref{Hyperbolic} hold. Then $\left( \bK _0\bftdiagup \bN
\right) \tilde{}=\{ 0\} $.
\end{Cor}

{\it Proof\/}: Suppose there are unitary functions $A_\pm \in \bL
^\infty (\bbR )$ such that $\overline{P_\pm \bK _0}=A_\pm \cH _\pm
$. Then this would violate the Schwarz-estimate (\ref{uJvBound}), and
therefore condition (PR3). Using irreducibility of $\p _+=\p _\hslash
$ and of $\p _-=\p _+\circ \t =\p _{-\hslash }$, we may reduce to
considering the cases when one of the spaces $\overline{P_\pm \bK _0}$
is $\bL ^2(\bbR )$. By Lemma \ref{Normal}, we are then back to the
case when $\bK _0$ or $\bK _0^{-1}$ is the graph of a densely defined
normal and dissipative operator $L$, or $L^{-1}$, respectively. We
will consider $L$ only. The other case goes the same way. Since
\begin{equation}\label{Pibf}
\left( \p _+\oplus \p _-\right) \left( 0,b,0\right) \left( f_+\oplus
f_-\right) (x)=e^{i\hslash bx}f_+(x)\oplus e^{-i\hslash bx}f_-(x)
\end{equation}
it follows that $L$ must anti-commute with the multiplication operator
$ix$ on $\bL ^2(\bbR )$. For deriving this, we used assumption (PR3)
at this point. We also showed in Lemma \ref{Normal} that $L$ must act
as a multiplication operator on the Fourier-transform side. But the
anti-commutativity is inconsistent with a known structure theorem in
\cite{Ped90}, specifically Corollary 3.3 in that paper. Hence there
are unitary functions $A_\pm $ in $\bL ^\infty (\bbR )$ such that
$\overline{P_\pm \bK _0}=A_\pm \cH _\pm$. But this possibility is
inconsistent with positivity in the form $\Re \ip{f_+}{f_-}\ge 0\,
,\quad \forall (f_+,f_-)\in \bK _0$ (see (\ref{JInnerProduct})) if
$\left( \bK _0\bftdiagup \bN \right) \tilde{}\ne \{ 0\} $. To see
this, note that $\bK _0$ is invariant under the unitary operators
(\ref{Pibf}) for $b\in \bbR _+$. The argument from Lemma
\ref{Heisenberg}, now applied to $\p _+\oplus \p _-$, shows that the
two subspaces
\[\bK _0^\infty :=\bigvee _{b\in \bbR }\left( \p _+\oplus \p _-\right)
\left( 0,b,0\right) \bK _0\]
and
\[\bK _0^{-\infty }:=\bigwedge _{b\in \bbR }\left( \p _+\oplus \p _-\right)
\left( 0,b,0\right) \bK _0\] are both invariant under the whole group
$\left( \p _+\oplus \p _-\right) (G)$. But the commutant of this is
$2$-dimensional: the only projections in the commutant are represented
as one of the following,
\[\left( \matrix{ 0 & 0\cr 0 & 0\cr }\right) ,\;\left( \matrix{
I & 0\cr 0 & 0\cr }\right) ,\;\left( \matrix{ 0 & 0\cr 0 & I\cr
}\right) ,\mathop{{\rm \ or}}\left( \matrix{ I & 0\cr 0 & I\cr
}\right) \, ,\] relative to the decomposition $\bL ^2(\bbR )\oplus \bL
^2(\bbR )$ of $\p _+\oplus \p _-$. The above analysis of the
anti-commutator rules out the cases $\left( \matrix{ I & 0\cr 0 & 0\cr
}\right) $ and $\left( \matrix{ 0 & 0\cr 0 & I\cr }\right) $, and if
$\left( \bK _0\bftdiagup \bN \right) \tilde{}\ne \{ 0\} $, we are left
with the cases $\bK _0^\infty =\{ 0\} $ and $\bK _0^\infty =\bL
^2(\bbR )\oplus \bL ^2(\bbR )$. Recall, generally $\bK _0^{-\infty
}\subset \bK _0\subset\bK _0^\infty $, as a starting point for the
analysis. A final application of the Beurling-Lax theorem (as in
\cite{LaPh}; see also \cite{DyMc70}) to (\ref{Pibf}) then shows that
there must be a pair of unitary functions $A_\pm $ in $\bL ^\infty
(\bbR )$ such that \begin{equation}\label{KAAHH} \bK _0=A_+\cH
_+\oplus A_-\cH _-
\end{equation}
where $\cH _\pm $ are the two Hardy spaces given by having $\hat{f}$
supported in $\left[ 0,\infty \right) $, respectively, $\left( -\infty
,0\right] $. The argument is now completed by noting that
(\ref{KAAHH}) is inconsistent with the positivity of $\bK _0$ in
(\ref{vJvPositivity}); that is, we clearly do not have $\ip{\left(
\matrix{ A_+h_+\cr A_-h_-\cr }\right) }{J\left( \matrix{ A_+h_+\cr
A_-h_-\cr }\right) }=2\Re \ip{A_+h_+}{A_-h_-}$ semidefinite, for all
$h_+\in \cH _+$ and all $h_-\in \cH _-$.  This concludes the proof of
the Corollary.\hfill$\Box$\medskip

\begin{Rem}
\label{UncorrelatedClosedSubspace}{\rm At the end of the above proof
of Corollary \ref{Positive}, we arrived at the conclusion
(\ref{KAAHH}) for the subspace $\bK _0$ under consideration. Motivated
by this, we define a closed subspace $\bK _0$ in a direct sum Hilbert
space $\bH _+\oplus \bH _-$ to be {\em uncorrelated} if there are
closed subspaces $\bD _\pm \subset \bH _\pm $ in the respective
summands such that
\begin{equation}
\bK _0=\bD _+\oplus \bD _-
\label{KDD}
\end{equation}
Contained in the corollary is then the assertion that every
semigroup-invariant $\bK _0$ in $\bL ^2(\bbR )\oplus \bL ^2(\bbR )$ is
uncorrelated, where the semigroup here is the subsemigroup $S$ in the
Heisenberg group $G$ given by
\begin{equation}
S=\left\{ (a,b,c)\mid b\in \bbR _+,\; a,c\in \bbR \right\} \, ,
\label{SinG}
\end{equation}
and the parameterization is the one from (\ref{NHeisenberg}). We also
had the representation $\p $ in the form $\p _+\oplus \p _-$ where the
respective summand representations $\p _\pm $ of $G$ are given by
(\ref{PiHBar}) relative to a pair $(\hslash ,-\hslash )$, $\hslash \in
\bbR \bftdiagdown \{ 0\} $ some fixed value of Planck's constant. In
particular, it is assumed in Corollary \ref{Positive} that each
representation $\p _\pm $ is {\em irreducible.} But for proving that
some given semigroup-invariant $\bK _0$ must be uncorrelated, this
last condition can be relaxed considerably; and this turns out to be
relevant for applications to Lax-Phillips scattering theory for the
wave equation with obstacle scattering \cite{LaPh}. In that context,
the spaces $\bD _\pm $ will be outgoing, respectively, incoming
subspaces; and the wave equation translates backwards, respectively
forwards, according to the unitary one-parameter groups $\p
_-(0,b,0)$, respectively, $\p _+(0,b,0)$, with $b\in \bbR $
representing the time-variable $t$ for the unitary time-evolution
one-parameter group which solves the wave equation under
consideration. The unitary-equivalence identity (\ref{JPiRelation})
stated before Lemma \ref{Operator} then implies equivalence of the
wave-dynamics before, and after, the obstacle scattering.  }
\end{Rem}

\npar Before stating our next result, we call attention to the
$(2n+1)$-dimensional Heisenberg group $G_n$ in the form $\bbR
^{2n+1}=\bbR ^n\times \bbR ^n\times \bbR $, in parameter form: $a,b\in
\bbR ^n$, $c\in \bbR $, and product rule
\[(a,b,c)\cdot (a^{\prime },b^{\prime },c^{\prime })=(a+a^{\prime
},b+b^{\prime },c+c^{\prime }+a\cdot b^{\prime })\] where $a+a^{\prime
}=(a_1+a_1^{\prime },\dots ,a_n+a_n^{\prime })$ and $a\cdot b^{\prime
}=\sum _{j=1}^na_jb_j^{\prime }$. For every (fixed) $b\in \bbR
^n\bftdiagdown \{ 0\} $, we then have a subsemigroup
\begin{equation}
S(b)=\left\{ (a,\b b,c)\mid \b \in \bbR _+,\; a\in \bbR ^n,\; c\in
\bbR \right\}\, ;
\label{SubsemigroupS}
\end{equation}
and we show in the next result that it is enough to have invariance
under such a semigroup in $G_n$, just for a single direction, defined
from some fixed $b\in \bbR ^n\bftdiagdown \{ 0\} $.

\begin{Satz}
\label{KNoughtUncorrelated}Let $\p _\pm $ be unitary representations
of the Heisenberg group $G$ on respective Hilbert spaces $\bH _\pm $,
and let $T:\bH _-\pf \bH _+$ be a unitary isomorphism which
intertwines $\p _-$ and $\p _+\circ \t $ as in {\rm (\ref{TPiMinus})}
where
\begin{equation}
\t (a,b,c)=(a,-b,-c)\, ,\quad \forall (a,b,c)\in G\simeq \bbR
^{2n+1}\, .
\label{TauInHeisenberg}
\end{equation}
Suppose there is $\hslash \in \bbR \bftdiagdown \{ 0\} $ such that 
\begin{equation}
\p _+(0,0,c)=e^{i\hslash c}I_{\bH _+}\, .
\label{PiPlusInHeisenberg}
\end{equation}
If $\bK _0\subset \bH _+\oplus \bH _-$ is a closed subspace which
is invariant under 
\[\left\{ (\p _+\oplus \p _-)(a,\b
b,c)\mid a\in \bbR ^n,\; \b \in \bbR _+,\; c\in \bbR \right\}\]
from
{\rm (\ref{SubsemigroupS}),} $b\in \bbR ^n\bftdiagdown \{ 0\} $, then we
conclude that $\bK _0$ must automatically be uncorrelated.
\end{Satz}

{\it Proof\/}: The group-law in the Heisenberg group yields the
following commutator rule:
\[(a,0,0)(0,b,0)(-a,0,0)=(0,b,a\cdot b)\]
for all $a,b\in \bbR ^n$. We now apply $\p =\p _+\oplus \p _-$ to
this, and evaluate on a general vector $f_+\oplus f_-\in \bK _0\subset
\bH _+\oplus \bH _-$: abbreviating $\p (a)$ for $\p (a,0,0)$, and $\p
(b)$ for $\p (0,b,0)$, we get
\[\p (a)\p (\b b)\p (-a)(f_+\oplus f_-)=e^{i\hslash \b a\cdot b}\p
_+(\b b)f_+\oplus e^{-i\hslash \b a\cdot b}\p _-(\b b)f_-\in \bK _0\]
valid for all $a\in \bbR ^n$, $\b \in \bbR _+$. Note, in
(\ref{PiPlusInHeisenberg}), we are assuming that $\p _+$ takes on some
specific value $e^{i\hslash c}$ on the one-dimensional center. Since
$\p _-$ is unitarily equivalent to $\p _+\circ \t $ by assumption (see
(\ref{PiPlusInHeisenberg})), we conclude that
\[\p _-(0,0,c)=e^{-i\hslash c}I_{\bH _-}\, ,\quad \forall c\in \bbR \, .\]
The argument really only needs that the two representations $\p _\pm $
define {\em different} characters on the center. (Clearly $\hslash \ne
-\hslash $ since $\hslash \ne 0$.) Multiplying through first with
$e^{-i\hslash \b a\cdot b}$, and integrating the resulting term
\[\p _+(\b b)f_+\oplus e^{-i2\hslash \b a\cdot b}\p _-(\b b)f_-\in \bK _0\]
in the $a$-variable, we get $\p _+(\b b)f_+\oplus 0\in \bK _0$. The
last conclusion is just using that $\bK _0$ is a closed subspace. But
we can do the same with the term
\[e^{i2\hslash \b a\cdot b}\p _+(\b b)f_+\oplus \p _-(\b b)f_-\in \bK _0\, ,\]
and we arrive at $0\oplus \p _-(\b b)f_-\in \bK _0$. Finally letting
$\b \pf 0_+$, and using strong continuity, we get $f_+\oplus 0$ and
$0\oplus f_-$ both in $\bK _0$. Recalling that $f_\pm $ are general
vectors in $P_\pm \bK _0$, we conclude that $P_+\bK _0\oplus P_-\bK
_0\subset \bK _0$, and therefore $\overline{P_+\bK _0}\oplus
\overline{P_-\bK _0}\subset \bK _0$. Since the converse inclusion is
obvious, we arrive at (\ref{KDD}) with $\bD _\pm =\overline{P_\pm \bK
_0}$.\hfill$\Box$\medskip

\npar The next result shows among other things that there are
representations $\p $ of the Heisenberg group $G_n$ (for each $n$)
such that the reflected representation $\p ^c$ of $G_n^c\simeq G_n$
(see Theorem \ref{PiCIrreducible}) acts on a nonzero Hilbert space
$\bH ^c=\left( \bK _0\bftdiagup \bN \right) \tilde{}$. However,
because of Lemma \ref{KernelPiC}, $\p ^c\left( G_n^c\right) $ will
automatically be an {\em abelian} group of operators on $\bH ^c$. To
see this, note that the proof of Theorem \ref{MDomain} shows that $\p
^c$ must act as the identity operator on $\bH ^c$ when restricted to
the one-dimensional center in $G_n^c\simeq G_n$.

\npar It will be convenient for us to read off this result from a more
general context: we shall consider a general Lie group $G$, and we fix
a right-invariant Haar measure on $G$.

\begin{Def}
\label{PositiveDefiniteDistribution}{\rm A distribution $F$ on the
Lie group $G$ will be said to be {\em positive definite} (PD) if
\newlength{\PDwidth} \settowidth{\PDwidth}{$\displaystyle \int _G\int
_GF(uv^{-1})\overline{f(u)}f(v)\,du\,dv\ge 0$} \newlength{\PDmargin}
\setlength{\PDmargin}{0.5\textwidth}
\addtolength{\PDmargin}{-0.5\PDwidth}
\[\int _G\int _GF(uv^{-1})\overline{f(u)}f(v)\,du\,dv\ge
0\mbox{\makebox[0pt][l]{\makebox[\PDmargin][r]{(PD)}\hss}}\] for all
$f\in C_c^\infty (G)$; and we say that $f$ is PD on some open subset
$\gO \subset G$ if this holds for all $f\in C_c^\infty (\gO )$. The
interpretation of the expression in (PD) is in the sense of
distributions. But presently measurable functions $F$ will serve as
the prime examples.

\npar We say that the distribution is {\em reflection-positive} (RP)
on $\gO $ ((RP$_\gO $) for emphasis) if, for some period-$2$
automorphism $\t $ of $G$, we have
\begin{equation}
F\circ \t =F
\label{FTau}
\end{equation}
and
\newlength{\RPwidth}
\settowidth{\RPwidth}{$\displaystyle \int _G\int _GF(\t
(u)v^{-1})\overline{f(u)}f(v)\,du\,dv\ge 0$}
\newlength{\RPmargin}
\setlength{\RPmargin}{0.5\textwidth}
\addtolength{\RPmargin}{-0.5\RPwidth}
\[\int _G\int _GF(\t (u)v^{-1})\overline{f(u)}f(v)\,du\,dv\ge
0\mbox{\makebox[0pt][l]{\makebox[\RPmargin][r]{(RP$_\gO $)}\hss}}\]
for all $f\in C_c^\infty (\gO )$.

\npar We say that some element $x$ in $G$ is (RP$_\gO $){\em
-contractive} if (RP$_\gO $) holds, and
\[0\le \int _G\int _GF(\t (u)v^{-1})\overline{f(ux)}f(vx)\, du\,dv\le
\int _G\int _GF(\t (u)v^{-1})\overline{f(u)}f(v)\,du\,dv\, ,\]
$\forall f\in C_c^\infty (\gO )$. Note that, since
\[\int _G\int _G\!F(\t (u)v^{-1})\overline{f(ux)}f(vx)\,du\,dv=\int
_G\int _G\!F(\t (u)\t (x)^{-1}xv^{-1})\overline{f(u)}f(v)\,du\,dv,\]
it follows that every $x$ in $H=G^\t $ is contractive: in fact,
isometric. If instead $\t (x)=x^{-1}$, then contractivity amounts to
the estimate
\[0\le \int _G\int _GF(\t (u)x^2v^{-1})\overline{f(u)}f(v)\,du\,dv\le
\int _G\int _GF(\t (u)v^{-1})\overline{f(u)}f(v)\,du\,dv\, ,\]
$\forall f\in C_c^\infty (\gO )$. Using the basic Lemma one can also
show that $x$ acts by contractions.}
\end{Def}

\npar The following result is useful, but an easy consequence of the
definitions and standard techniques for positive definite
distributions; see for example \cite{Jor88,Pra89}.

\begin{Satz}
\label{PiExtends}Let $F$ be a distribution on a Lie group $G$ with
a period-$2$ automorphism $\t $, and suppose $F$ is $\t $-invariant,
{\rm (PD)} holds on $G$, and {\rm (RP$_\gO $)} holds on some open, and
semigroup-invariant, subset $\gO $ in $G$. Then define
\[(\p (u)f)(v):=f(vu)\, ,\quad \forall u,v\in G,\; \forall f\in
C_c^\infty (G)\, ;\]
and
\[Jf:=f\circ \t \, .\]
Let $\bH (F)$ be the Hilbert space obtained from the GNS construction,
applied to {\rm (PD),} with inner product on $C_c^\infty (G)$ given by
\[\ip{f}{g}:=\int _G\int _GF(uv^{-1})\overline{f(u)}g(v)\,du\,dv\, .\]
Then $\p $ extends to a unitary representation of $G$ on $\bH (F)$,
and $J$ to a unitary operator, such that
\[J\p =(\p \circ \t )J\, .\]
If {\rm (RP$_\gO $)} further holds, as described, then $\p $ induces
{\rm (}via Theorem {\rm \ref{PiCIrreducible})} a unitary
representation $\p ^c$ of $G^c$ acting on the new Hilbert space $\bH
^c$ obtained from completing in the new inner product from {\rm
(RP$_\gO $),} and dividing out with the corresponding kernel.
\end{Satz}

\npar The simplest example of a function $F$ on the Heisenberg group
$G_n$ satisfying (PD), but not (RP$_\gO $), for nontrivial $\gO $'s,
may be obtained from the Green's function for the sub-Laplacian on
$G_n$; see \cite[p. 599]{Ste93} for details.

\npar If complex coordinates are introduced in $G_n$, the formula for
$F$ takes the following simple form: let $z\in \bbC ^n$, $c\in \bbR $,
and define
\[F(z,c)=\frac 1{\left( \left| z\right| ^4+c^2\right) ^n}\, .\]
Then we adapt the product in $G_n$ to the modified definition as
follows:
\[(z,c)\cdot (z^{\prime },c^{\prime })=(z+z^{\prime },c+c^{\prime
}+\left\langle z,z^{\prime }\right\rangle )\quad \forall z,z^{\prime
}\in \bbC ^n,\; \forall c,c^{\prime }\in \bbR \, ,\] where
$\left\langle z,z^{\prime }\right\rangle $ is the symplectic form
\[\left\langle z,z^{\prime }\right\rangle :=2\Im (z\cdot
\halo{z}^{\prime })\, .\]
The period-$2$ automorphism $\t $ on $G_n$ we take as 
\[\t (z,c)=(\halo{z},-c)\]
with $\halo{z}$ denoting complex conjugation, i.e., if $z=(z_1,\dots
,z_n)$, then $\halo{z}:=(\halo{z}_1,\dots ,\halo{z}_n)$.

\npar The simplest example where both (PD) and (RP$_\gO $) hold on the
Heisenberg group $G_n$ is the following:

\begin{Ex}
\label{SubLaplacian}{\rm Let $\z =(\z _1,\dots \z _n)\in \bbC ^n$,
$\x _j=\Re \z _j$, $\h _j=\Im \z _j$, $j=1,\dots ,n$. Define
\[F(z,c)=\int _{\bbR ^{2n}}\frac {e^{i\Re (z\cdot \halo{\z })}}{\prod
_{j=1}^n(\left| \z _j\right| ^2+1)}\,d\x _1\cdots d\x _n\,d\h _1\cdots
d\h _n\, .\] Let $\gO :=\{ (z,c)\in G_n\mid z=(z_j)_{j=1}^n,\; \Im
z_j>0\} $.  Then (PD) holds on $G_n$, and (RP$_\gO $) holds, referring
to this $\gO $. Since the expression for $F(z,c)$ factors, the problem
reduces to the $(n=1)$ special case. There we have
\[F(z,c)=\int _{\bbR ^2}\frac {e^{i(x\x +y\h )}}{\x ^2+\h ^2+1}\,d\x
\,d\h \, ;\] and if $f\in C_c^\infty (\gO)$ with $\gO =\left\{
(z,c)\mid y>0\right\} $, then
\begin{eqnarray*}
&&\int _{G_1}\int _{G_1}F(\t(u)v^{-1})\overline{f(u)}f(v)\,du\,dv\\
&&=\int _{\bbR ^8}\frac{e^{i(x-x^{\prime })\x }\,e^{-i(y+y^{\prime
})\h }}{\x ^2+\h ^2+1}\overline{f(x+iy,c)}f(x^{\prime }+iy^{\prime
},c^{\prime })\,d\x \,d\h \,dx\,dy\,dc\,dx^{\prime }\,dy^{\prime
}\,dc^{\prime }\, .
\end{eqnarray*}
Let $\tilde{f}$ denote the Fourier transform in the $x$-variable,
keeping the last two variables $(y,c)$ separate. Then the integral
transforms as follows:
\[\int _{\bbR ^5}\frac {e^{-(y+y^{\prime })\sqrt{1+\x ^2}}}{\sqrt{1+\x
^2}}\overline{\tilde{f}(\x ,y,c)}\tilde{f}(\x ,y^{\prime },c^{\prime
})\,d\x \,dy\,dy^{\prime }\,dc\,dc^{\prime }\, .\] Introducing the
Laplace transform in the middle variable $y$, we then get (since $f$
is supported in $y>0$)
\[\int _0^\infty e^{-y\sqrt{1+\x ^2}}\tilde{f}(\x ,y,c)\,dy=\tilde{f}_\l
(\x ,\sqrt{1+\x ^2},c)\, ;\] the combined integral reduces further:
\[\int _\bbR \left| \vphantom{\int \tilde{f}_\l (\x ,\sqrt{1+\x
^2},c)\,dc}\smash{\int _\bbR \tilde{f}_\l (\x ,\sqrt{1+\x
^2},c)\,dc}\right| ^2\frac {d\x}{1+\x ^2}\] which is clearly positive;
and we have demonstrated that (RP$_\gO $) holds. It is immediate that
$F$ is $\t $-invariant (see (\ref{FTau})), and also that it satisfies
(PD) on $G_n$.  }
\end{Ex}

\section{\protect\label{axb}The $(ax+b)$-Group}
\setcounter{equation}{0}

We showed that in general we get a unitary representation $\p^c$ of
the group $G^c$ from an old one $\p$ of $G$, provided $\p$ satisfies
the assumptions of reflection positivity.  The construction as we saw
uses a certain cone $C$ and a semigroup $H\exp C$, which are part of
the axiom system. What results is a new class of unitary
representations $\p^c$ satisfying a certain spectrum condition
(semi-bounded spectrum).

\npar But, for the simplest non-trivial group $G$, this semi-boundness
turns out {\it not} to be satisfied in the general case.  Nonetheless,
we still have a reflection construction getting us from unitary
representations $\p$ of the $(ax+b)$-group, such that $\p \circ \t
\simeq \p$ (unitary equivalence), to associated unitary
representations $\p^c$ of the same group. The (up to conjugation)
unique non-trivial period-$2$ automorphism $\t$ of $G$, where $G$ is
the $(ax+b)$-group, is given by
\[\t(a,b) = (a,-b)\, . \] 
Recall that the $G$ may be identified with the matrix-group
\[\left\{\left. \left(\matrix{a & b\cr 0 & 1\cr}\right)\,\right|
\,a>0, b\in \bbR\right\}\]
and $(a,b)$ corresponds to the matrix
$ \left(\matrix{a & b\cr 0 & 1\cr}\right)$.
In this realization the Lie algebra of $G$ has
the basis
\[X = \left(\matrix{1 & 0\cr 0 & 0\cr}\right)
\quad \mbox{and}\quad Y= \left(\matrix{0 & 1\cr 0 & 0 
\cr}\right)\, .\]
We have  $\exp (tX) = (e^t,0)$ and $\exp (sY)
= (1,s)$. Hence $\t (X) = X$ and
$\t (Y) = -Y$. Thus
$\fh = \bbR X$ and
$\fq = \bbR Y$.
We notice the commutator relation $[X,Y] = Y$.
The possible $H$-invariant
cones in $\fq$ are
$\pm \{tY\mid t\ge 0\}$.
It is known
from Mackey's theory that $G$ has two inequivalent, unitary,
irreducible, infinite-dimensional representations $\p_{\pm}$, 
and
it is immediate that we have the unitary equivalence
(see details below):
\begin{equation}\label{E:5.1}
\p_+\circ\t \simeq \p_-\, .
\end{equation}
Hence, if we set $\p := \p_+\oplus \p_-$, then
$\p\circ \t \simeq \p$, so we have the setup
for the general theory. We show that $\p$ may be
realized on $\bL^2(\bbR)\oplus \bL^2(\bbR)
\simeq \bL^2(\bbR,\bbC^2)$, and we find and
classify the invariant positive subspaces
$\bK_0\subset \bL^2(\bbR,\bbC^2)$. To understand
the interesting cases for the $(ax+b)$-group $G$, we need
to relax the invariance condition: We shall {\it not}
assume invariance of $\bK_0$ under
the semigroup $\{\p (1,b)\mid b\ge 0\}$, but only
under the 
infinitesimal unbounded generator $\p (Y)$. With this,
we still get the correspondence $\p \mapsto \p^c_{\bK_0}$ as
described above.

\npar
We use the above notation. We know from Mackey's theory
\cite{Mac} that there are two inequivalent irreducible
infinite-dimensional representations of $G$, and we
shall need them in the following alternative formulations: Let
$\cL_{\pm}$ denote the respective Hilbert space
$\bL^2(\bbR_{\pm})$ with
the multiplicative invariant measure $d\m_{\pm} =
dp/|p|$, $p\in \bbR_\pm$. Then the formula
\begin{equation}\label{E:5.4}
f\mapsto e^{ipb}f(pa)
\end{equation}
for functions $f$ on $\bbR$ restricts to two unitary irreducible
representations, denoted by $\p_\pm$ of $G$ on
the respective spaces $\cL_\pm$. Let
$Q(f) (p) := f(-p)$ denote the canonical mapping
from $\cL_+$ to $\cL_-$, or equivalently
from $\cL_-$ to $\cL_+$. Then we have
for $g\in G$ (cf. (\ref{E:5.1})):
\begin{equation}\label{E:5.5}
Q\p_+(g) = \p_-(\t (g))Q
\end{equation}
For the representation $\p := \p_+\oplus \p_-$ on
$\bH := \cL_+\oplus \cL_-$ we therefore have
\begin{equation}\label{E:5.6}
J\p (g) = \p(\t (g))J,\quad g\in G
\, ,
\end{equation}
where $J$ is the unitary 
involutive operator on $\bH$ given by
\begin{equation}\label{E:5.7}
J = \left(\matrix{ 0 & Q\cr Q & 0\cr}\right)\, .
\end{equation}
Instead of the above $p$-realization of $\p$ we
will mainly use the following $x$-formalism. The
map $t\mapsto \pm e^t$ defines an isomorphism
$L_\pm : \cL_\pm \pf \bL^2(\bbR)$, where we use the
(additive) Lebesgue measure $dx$ on $\bbR$.
For $g = (e^s,b)\in G$ and $f\in \bL^2(\bbR)$, set
\begin{equation}\label{E:5.8}
(\p_\pm (g)f)(x) := e^{\pm ie^xb}f(x+s), \quad x\in \bbR\, .
\end{equation}
A simple calculation shows that $L_\pm$ intertwines the
old and new construction of $\p_\pm$, excusing our
abuse of notation. In this realization $Q$ becomes simply
the identity operator $Q(f)(x) = f(x)$.
The involution $J : \bL^2(\bbR,\bbC^2)$ is now simply given by
\[J(f_0,f_1) = (f_1,f_0)\]
or $J = \left(\matrix{0 & 1\cr 1 & 0\cr}\right)$.

\npar
In this formulation the operator
\begin{equation}\label{E:5.9}
L := \p_\pm(\D_H - \D_q) = \p_\pm(X^2-Y^2)
\end{equation}
takes the form
\begin{equation}\label{E:5.10}
L = \left(\frac{d\, }{dx}\right)^2\, + \, e^{2x}\, ,
\end{equation}
but it is on $\bL^2(\bbR)$ and
$-\infty < x<\infty$. This operator is known to
have defect indices $(1,0)$
\cite{Jor75,NeSt59}, which means that it cannot be
extended to a selfadjoint operator on 
$\bL^2(\bbR)$. Using a theorem from
\cite{Jor75,ReSi75} we can see this by comparing the
quantum mechanical problem for a particle
governed
by $-L$ as a
Schr\"odinger operator (i.e., a strongly repulsive force) 
with the corresponding classical
one governed (on each energy surface) by
\[ E_{\rm kin}+E_{\rm pot}=\left(\frac{dx}{dt}\right)^2\, - \, e^{2x} = E\, .\]
The escape time for this particle to
$x = \pm \infty$ is 
\begin{equation}\label{E:5.11}
t_{\pm} = \int_{\rm finite}^{\pm 
\infty}\frac{dx}{\sqrt{E+e^{2x}}\, }\, ,
\end{equation}
i.e., $t_\infty$ is finite, and
$t_{-\infty} = \infty$. We elaborate on this point below.
The nonzero defect vector for
the quantum mechanical problem
corresponds to  a
boundary condition at $x=\infty$ since this is
the singularity which is reached in finite time.

\npar The fact from \cite{Jor75} we use for the defect index assertion
is this: The Schr\"{o}dinger operator $H=-\left(\frac{d\,
}{dx}\right)^2\, + V(x)$ for a single particle has nonzero defect
solutions $f_\pm \in \bL ^2(\bbR )$ to $H^*f_\pm =\pm if_\pm $ iff
there are solutions $t\mapsto x(t)$ to the corresponding classical
problem
\[ E=\left(\frac{dx(t)}{dt}\right ) ^2+V(x(t))\]
with finite travel-time to $x=+\infty $, respectively, $x=-\infty
$. The respective (possibly infinite) travel-times are
\[t_{\pm \infty } = \int_{\rm finite}^{\pm 
\infty }\frac{dx}{\sqrt{E-V(x)}}\, .\] The correspondence principle
states that one finite travel-time to $+\infty $ (say) yields a
dimension in the associated defect space, and similarly for the other
travel-time to $-\infty $.

\npar In the $x$-formalism, (\ref{E:5.5}) from above then simplifies
to the following identity for operators on the {\it same} Hilbert
space $\bL^2(\bbR)$ (carrying the two inequivalent representations
$\p_+$ and $\p_-$):
\begin{equation}\label{E:5.12}
\p_+(g) = \p_-(\t (g)),\quad g\in G\, .
\end{equation}
We realize the representation $\p = \p_+\oplus \p_-$ in
the Hilbert space
$\bH = \bL^2(\bbR)\oplus \bL^2(\bbR) = \bL^2(X_2)
$
where $X_2 = {0}\times \bbR\cup {1}\times \bbR$. We
may represent $J$ by an automorphism $\th : X_2\pf X_2$
(as illustrated in Proposition \ref{P:3.3}):
\[
\th (0,x) := (1,x)\quad\mbox{and}\quad \th (1,y) = (0,y)\, ,\quad
x,y\in \bbR\, ,\]
and
\[J(f)(\go ) = f(\th (\go))\, ,\quad \go \in X_2\, .\]
Notice that the subset
\[X_2^\th = \{\go \in X_2\mid \th (\go) = \go\}\]
is empty. Define for $f\in \bL^2(X_2)$, $f_k(x) = f(k,x)$,
$k=0,1$, $x\in \bbR$. We have for
$g = (e^s,b)\in G$:
\[ \left(\p (g)f\right)_0(x) = e^{ibe^x}f_0(x+s)=(\halo{\p }_+(g)f_0)(x)\]
and
\[ \left(\p (g)f\right)_1(x) = e^{-ibe^x}f_1(x+s)=(\halo{\p }_-(g)f_1)(x)\, .\]

\begin{Prop}\label{Negative}Let $\p =\p _+\oplus \p _-$ be the
representation from {\rm (\ref{E:5.1})--(\ref{E:5.6})} above of the
$(ax+b)$-group $G$. Then the only choices of reflections $\bK_0$ as in
Remark \ref{WhatToAssume} for the sub-semigroup
$S=\left\{\left(a,b\right)\in G\mid b>0\right\}$ will have
$\bK=\left(\bK_0\bftdiagup \bN\right)\tilde{}$ equal to $0$.
\end{Prop}

{\it Proof\/}: Let $\bK_0$ be as specified in Remark
\ref{WhatToAssume} relative to the semigroup $S$, and let $P_{\bK_0}$
be the representation of the corresponding orthogonal projection
operator as given in (\ref{E:5.14})--(\ref{QJQPositivity}) in terms of
the measurable field $\bbR\ni \x \mapsto Q(\x ) $. Specifically, the
space $\bK_0\subset \bH $ (with the positivity and invariance
properties from Section 2) will then be translation invariant, i.e.,
invariant under the translation group
\begin{equation}\label{E:5.13}
\left(\matrix{f_0(x)\cr f_1(x)\cr}\right)\mapsto
\left(\matrix{f_0(x+s)\cr f_1(x+s)\cr}\right)\, ,
\quad x,y,s\in \bbR\, .
\end{equation}
Hence the projection in $\bL^2(\bbR,\bbC^2)$ onto
$\bK_0$, denoted by $P_{\bK _0}$, may be represented as a
multiplication operator in the Fourier transform space
\[f = \left(\matrix{f_0\cr f_1\cr}\right)\, ,
\quad
\hat{f}(\x ) =\left(\matrix{\hat{f_0}(\x )\cr
 \hat{f_1}(\x )\cr}\right)\, ,\quad \x \in \bbR\, ,\]
where as usual
\[\hat{f_k}(\x ) = \frac{1}{\sqrt{2\pi}}\int_{\bbR}e^{-i\x 
x}f_k(x)\,dx\, ,\quad k=0,1\, .\]

\begin{Lemma}\label{QField}Let $Q$ be the projection in $\bL
^2(\bbR )\oplus \bL ^2(\bbR )$ onto a translation-invariant
$J$-positive subspace. Then $Q$ is represented by a measurable field
of $2\times 2$ complex matrices $\bbR \ni \x \mapsto \left( Q_{ij}(\x
)\right) _{ij=1}^2$ such that $\left|Q_{12}(\x )\right| ^2=Q_{11}(\x
)Q_{22}(\x )$ a.e. on $\bbR$, and $Q_{12}(\x )+Q_{21}(\x )\ge 0$ a.e.;
and conversely.
\end{Lemma}

{\it Proof\/}: Since all the operators commuting with the translation
group (\ref{E:5.14}) are known (see, e.g., \cite{LaPh}), there is a
measurable field of projections $Q (\x ) : \bbC^2\pf \bbC^2$, i.e.,
$Q(\x )^2 = Q(\x ) = Q(\x )^*$, such that
(\ref{E:5.14})--(\ref{QJQPositivity}) hold:
\begin{equation}\label{E:5.14}
\left( P_{\bK_0}f\right) ^\wedge (\x ) = Q(\x )\hat{f}(\x )\, .
\end{equation}
With $J = \left(\matrix{ 0 & 1\cr 1 & 0\cr}\right)
: \bbC^2\pf \bbC^2$ as
before, we have the basic positivity:
\begin{equation}\label{QJQPositivity}
Q(\x )J Q(\x ) \ge 0\, ,\quad 
\x \in \bbR\, .
\end{equation}
Hence
\newlength{\trqjqgez}
\settowidth{\trqjqgez}{$\displaystyle \Tr (Q(\x )J Q(\x )) \ge 0$}
\newlength{\bmargin}
\setlength{\bmargin}{0.5\textwidth}
\addtolength{\bmargin}{-0.5\trqjqgez}
\[\phantom{\Tr (Q(\x )J Q(\x ))}\mbox{\makebox[0pt][r]{\hss$\displaystyle
\det(Q(\x )JQ(\x ))$}}\ge
0\mbox{\makebox[0pt][l]{\makebox[\bmargin][r]{(a)}\hss}}\]
and
\[\Tr (Q(\x )J Q(\x )) \ge 0\mbox{\makebox[0pt]{\,
.\hss}\makebox[0pt][l]{\makebox[\bmargin][r]{(b)}\hss}}\]

\npar Since $\det(QJQ)=-\det(Q)=-\det (Q^2)=-(\det Q^2)\le 0$, it
follows from (a) that $\det Q (\x ) = 0$, and, from (a)--(b), that
$Q(\x )$ is for each $\x $ a projection into a subspace in $\bbC^2$ of
dimension $0$ or $1$. Write $Q=(Q_{ij})$, with $Q_{ij}:\bbR \pf \bbC $
measurable. Then $Q=Q^*$ gives, for $\x \in \bbR $,
\[Q_{11}(\x ),Q_{22}(\x )\in \bbR \und Q_{21}(\x )=\overline{Q_{12}(\x )}\, .\]
The relation $Q^2(\x )=Q(\x )$ implies
\[Q_{11}(\x )^2+\left| Q_{12}(\x )\right| ^2=Q_{11}(\x )\, ,\]
\[Q_{22}(\x )^2+\left| Q_{12}(\x )\right| ^2=Q_{22}(\x )\, ,\]
and
\[\left( Q_{11}(\x )+Q_{22}(\x )\right) Q_{12}(\x )=Q_{12}(\x )\, .\]
In particular
\[0\le Q_{11}(\x ),Q_{22}(\x )\le 1\]
and
\[\left| Q_{12}(\x )\right| ^2=Q_{11}(\x )\left( 1-Q_{11}(\x )\right)
=Q_{22}(\x )\left( 1-Q_{22}(\x )\right) \, .\]
{F}rom $\det Q(\x )=0$, we finally get
\[\left| Q_{12}(\x )\right| ^2=Q_{11}(\x )Q_{22}(\x )\, .\]

\begin{Cor}\label{QMatrix}
These relations imply the following for the matrix $Q$:
\begin{enumerate}
\item[\hss\llap{\rm 1)}] If $Q_{12}(\x )=0$ then we have the three
possibilities:
\begin{eqnarray*}
Q(\x ) &=& 0\, ,\\
Q(\x ) &=& \left( \matrix{ 1 & 0\cr 0 & 0\cr}\right) \, ,\mbox{ and }\\
Q(\x ) &=& \left( \matrix{ 0 & 0\cr 0 & 1\cr}\right) \, .
\end{eqnarray*}
In all those cases, we have $Q(\x )JQ(\x )=0$.
\item[\hss\llap{\rm 2)}] If $Q_{12}(\x )\ne 0$, then $0<Q_{22}(\x
)=1-Q_{11}(\x )<1$. Let $\m (\x )=Q_{12}(\x )/Q_{11}(\x )$. Then by
$\Tr (Q(\x )J Q(\x )) \ge 0$ we have $\Re \m (\x )\ge 0$ and
\begin{equation}\label{QMatrixFormula}
Q(\x )=\frac 1{1+\left| \m (\x )\right| ^2}\left( \matrix{ 1 & \m (\x
)\cr \overline{\m (\x )} & \left| \m (\x )\right| ^2\cr }\right) \, .
\end{equation}
With $\l =\bar{\m }$ we get that the image of $Q(\x )$ is given by
\[\left\{ \left. u(\x )\left( \matrix{ 1\cr \l (\x )\cr }\right)
\,\right| \,u(\x )\in \bbC \right\} \, .\] Specifying to our
situation, $f=\left( \matrix{ f_0\cr f_1\cr }\right) \in \bK_0$ if and
only if
\begin{equation}\label{flamf}
\hat{f}_1(\x )=\l (\x )\hat{f}_0(\x )\, .
\end{equation}
\end{enumerate}
\end{Cor}

Since $Q(\x )$ is a measurable field of projections, the function
$\bbR \ni \x \mapsto \l (\x )$ must be measurable, but it may be
unbounded. This also means that $P_{\bK_0}$ is the projection onto the
graph of the operator $T_0:f_0\mapsto f_1$ where $f_0$ and $f_1$ are
related as in (\ref{flamf}), and the Fourier transform $\hat{\cdot }$
is in the $\bL^2$-sense.
\smallskip

\npar {\it Proof of Proposition \ref{Negative} {\rm (}continued\/{\rm
)}} : We first assume that $\bK_0$ arises this way as the graph of an
operator $T_0$ as described. This assumption will then be ``removed''
later.

\npar The assumed invariance of $\bK_0$ under $\p =\p _+\oplus \p _-$
takes the form
\begin{equation}\label{InvarKNought}
\left( \matrix{ \p _+(b) & 0\cr 0 & \p _-(b)\cr }\right) \bK_0\subset
\bK_0\, ,\quad \forall b>0\, .
\end{equation}
Let $\cD \subset \bL ^2(\bbR )$ consist of the $\bL ^2(\bbR )$-closure
of the functions $f_0$ such that $\ip{\left( \matrix{ \hat{f}_0\cr \l
\hat{f}_0\cr }\right) }{J\left( \matrix{ \hat{f}_0\cr \l \hat{f}_0\cr
}\right) }=0$. This may also be expressed in the form
\begin{equation}\label{NullFHat}
\int _{-\infty }^\infty \Re \l (\x )\left| \hat{f}_0(\x )\right|
^2\,d\x =0\, .
\end{equation}
It follows from (\ref{InvarKNought}) and Lemma \ref{BasicLemma} (the
Basic Lemma) that
\[\g _b:=\left( \matrix{ \p _+(b) & 0\cr 0 & \p _-(b)\cr }\right)
=\left( \matrix{ \p _+(b) & 0\cr 0 & \p _+(-b)\cr }\right) \]
for $b>0$ satisfies
\begin{equation}\label{GammaBound}
\ip{\g _b(v)}{J\g _b(v)}\le \ip{v}{Jv}
\end{equation}
for all $v\in \bK_0$ and $b\in \bbR _+$. When the explicit operators
are substituted into the latter estimate, we get
\[\int _{-\infty }^\infty \Re \l (\x ) \left| \left( \p _+(b)f_0\right)
^\wedge (\x )\right| ^2\,d\x \le \int _{-\infty }^\infty \Re \l
(\x )\left| \hat{f}_0(\x )\right| ^2\,d\x \, ,\]
valid for $b\in \bbR _+$, and
\[\left( \matrix{ f_0\cr T_0f_0\cr }\right) \in \bK_0\subset\left(
\matrix{ \bL ^2(\bbR )\cr \bL ^2(\bbR )}\right) \, .\] It follows that
$\p _+(b)$ maps the subspace $\cD $ into itself when $b\in \bbR _+$;
and, as a consequence, the Lax-Phillips setup applies to $\cD $ as a
closed subspace in $\bL ^2(\bbR )$, relative to the unitary
one-parameter group $\left\{ \p _+(b)\mid b\in \bbR \right\} $ of
operators in $\bL ^2(\bbR )$. Let
\begin{eqnarray}
\cD _\infty &:=&\bigvee _{b\in \bbR }\p _+(b)\cD\, ,\label{DInfinitybis}\\
\cD _{-\infty } &:=&\bigwedge _{b\in \bbR }\p _+(b)\cD \,
,\label{DMinusInfinitybis}
\end{eqnarray}
where $\bigvee $ and $\bigwedge $ denote the lattice operations on
closed subspaces in $\bL ^2(\bbR )$, and
\[\left( \p _+(b)f\right) (x)=e^{ibe^x}f(x)\, ,\quad f\in \bL
^2(\bbR ),\; b,x\in \bbR .\] It follows from the ansatz
(\ref{DInfinitybis})--(\ref{DMinusInfinitybis}) that both of the
spaces $\cD _\infty $ and $\cD _{-\infty }$ are invariant under
$\left\{ \p _+(b)\mid b\in \bbR \right\}$, and moreover that
\begin{equation}\label{DNest}
\cD _{-\infty }\subset \cD \subset \cD _\infty \, .
\end{equation}
It is enough to show that the assumption $\cD \ne \{ 0\} $ leads to a
trivial quotient space $\left( \bK_0\bftdiagup \bN \right)
\tilde{}$. Let
\[\t (s)f(x)=f(x+s)\, ,\quad f\in \bL ^2(\bbR ),\; s,x\in \bbR \]
be the translation part. We have
\begin{equation}\label{TauPi}
\t (s)\p _+(b)=\p _+(e^sb)\t (s)
\end{equation}
and we conclude that $\cD _{\pm \infty }$ are also both invariant
under $\left\{ \t (s)\mid s\in \bbR \right\} $. Since, as we noted,
the system (\ref{TauPi}) is irreducible in $\bL ^2(\bbR )$, we
conclude by Schur's lemma that $\cD _\infty =\bL ^2(\bbR )$. Recall
$\cD \ne 0$ was assumed at the outset. For the space $\cD _{-\infty
}$, we then have only two possibilities, $\cD _{-\infty }=\{ 0\} $ and
$\cD _{-\infty }=\bL ^2(\bbR )$, again by Schur's lemma, and the first
possibility must be ruled out by virtue of the Lax-Phillips theorem
\cite{LaPh}. Notice that the spectrum of $\left\{ \p _+(b)\mid b\in
\bbR \right\} $ is evidently a half-line, and the two properties, $\cD
_{-\infty }=\{ 0\} $ and $\cD _\infty =\bL ^2(\bbR )$, would
contradict the conclusion in the Lax-Phillips theorem, to the effect
that the spectrum would then necessarily have to be two-sided, i.e.,
all $\bbR =(-\infty ,\infty )$, and of homogeneous Lebesgue type,
i.e., unitarily equivalent, up to multiplicity, with translation on
the line.

\npar Only the possibility $\cD _{-\infty }=\bL ^2(\bbR )$ remains to
be considered. But we have
\[\cD _{-\infty }\subset \cD \subset P_0\bK _0\, ,\]
so it would follow that $\cD =\bL ^2(\bbR )$, and we are then reduced
back again to the case $Q=P_{\bK _0}=\left( \matrix{ 1 & 0\cr 0 & 0\cr
}\right) $ from part I of the present proof; i.e., to a trivial
induced Hilbert space $\left( \bK _0\bftdiagup \bN \right) \tilde{}$
as already noted.\hfill$\Box$\medskip

\npar The following argument deals with the general case, avoiding the
separation of the proof into the two cases (I) and (II): If vectors
$v\in \bK _0$ are expanded as $v=\left( \matrix{ h\cr k\cr }\right) $,
$h=Q_{11}h+Q_{12}k$, $k=Q_{21}h+Q_{22}k$, we can introduce $\cD
=\left\{ h\in \bL ^2(\bbR )\,\left| \, \exists k\in \bL ^2(\bbR )
\mathop{{\rm s.t.}}\left( \matrix{ h\cr k\cr }\right) \in \bN
\right. \right\} $. If $b>0$, we then have from (\ref{InvarKNought}):
\begin{eqnarray*}
\p _+(b)h &=& Q_{11}\p _+(b)h+Q_{12}\p _+(-b)k\, ,\\
\p _+(-b)k&=&Q_{21}\p _+(b)h+Q_{22}\p _+(-b)k\, ,
\end{eqnarray*}
valid for any $\left( \matrix{ h\cr k\cr }\right) \in \bK _0$, and
$b\in \bbR _+$. So it follows from Lemma \ref{BasicLemma} again that
$\cD $ is invariant under $\left\{ \p _+(b)\mid b>0\right\} $, and
also under the whole semigroup $\left\{ \p _+(g)\mid g\in S\right\} $
where $\p _+$ is now denoting the corresponding (see (\ref{TauPi}))
irreducible representation of $G$ on $\bL ^2(\bbR )$. Hence, we may
apply the Lax-Phillips argument to the induced spaces $\cD _{\pm
\infty }$ from (\ref{DInfinitybis})--(\ref{DMinusInfinitybis}). If
$\left( \bK _0\bftdiagup \bN \right) \tilde{}$ should be $\ne \{ 0\}
$, then $\cD =\{ 0\} $ by the argument. Since we are assuming $\left(
\bK _0\bftdiagup \bN \right) \tilde{}\ne\{ 0\} $, we get $\cD =\{ 0\}
$, and as a consequence the following operator graph representation
for $\bK _0$: $\left( \bK _0\bftdiagup \bN\right) \tilde{}=\b \left(
G(L)\right) $ where $G(L)$ is the graph of a closed operator $L$ in
$\bL ^2(\bbR )$.  Specifically, this means that the linear mapping
$\bK _0\bftdiagup \bN \ni \left( \matrix{ h\cr k\cr }\right) +\bN
\mapsto h$ is well-defined as a linear closed operator. This in turn
means that $\bK _0$ may be represented as the graph of a closable
operator in $\bL ^2(\bbR )$ as discussed in the first part of the
proof. Hence such a representation could have been assumed at the
outset.

\begin{Rem}\label{BorchersCMP92}{\rm 
In a recent paper on local quantum field theory \cite{Bor92}, Borchers
considers in his Theorem II.9 a representation $\p $ of the
$(ax+b)$-group $G$ on a Hilbert space $\bH $ such that there is a {\em
conjugate linear} $J$ (i.e., a period-$2$ antiunitary) such that $J\p
J=\p \circ \t $ where $\t $ is the period-$2$ automorphism of $G$
given by $\t (a,b):=(a,-b)$. In Borchers's example, the one-parameter
subgroup $b\mapsto \p (1,b)$ has semibounded spectrum, and there is a
unit-vector $v_0\in \bH $ such that $\p (1,b)v_0=v_0$, $\forall b\in
\bbR $. The vector $v_0$ is cyclic and separating for a von Neumann
algebra $M$ such that $\p (1,b)M\p (1,-b)\subset M$, $\forall b\in
\bbR _+$. Let $a=e^t$, $t\in \bbR $. Then, in Borchers's construction,
the other one-parameter subgroup $t\mapsto \p (e^t,0)$ is the modular
group $\D ^{it}$ associated with the cyclic and separating vector
$v_0$ (from Tomita-Takesaki theory \cite[vol. I]{BrRo}). Finally, $J$
is the corresponding modular conjugation satisfying $JMJ=M^{\prime }$
when $M^{\prime }$ is the commutant of $M$.  }
\end{Rem}
\medskip
{\bf Acknowledgements. } The authors would like to thank the
University of Iowa (May 1994), the Louisiana State University (Nov.\
1995), the University of Odense (June 1995), and the Mittag-Leffler
Institute (Spring 1996) for their support. We would also want to thank
Bent {\O}rsted and Steen Pedersen for helpful discussions. We would
like to thank Brian Treadway for his excellent typesetting of the
final version of the paper in \LaTeX . Both authors were also
supported in part by the National Science Foundation.

\end{document}